\pdfoutput=1
\documentclass[reqno,a4paper,11pt]{amsart}

\usepackage[T1]{fontenc}
\usepackage{lmodern}
\usepackage{bbm}
\usepackage{microtype}
\usepackage{mathtools,amssymb}
\usepackage{xcolor}
\usepackage{graphicx}
\usepackage{enumitem}
\usepackage[margin=1in]{geometry}

\usepackage[pdfusetitle,colorlinks=true,linkcolor=blue!45!black,citecolor=blue!45!black,urlcolor=blue!45!black]{hyperref}

\usepackage{tikz}
\usetikzlibrary{patterns,arrows.meta,patterns.meta,calc,positioning,shadings}
\tikzset{
 paperfigure/.style={font=\small, line cap=round, line join=round},
 figaxis/.style={draw=black!65, line width=.4pt,
   -{Stealth[length=1.4mm,width=1.05mm]}},
 figoutline/.style={draw=black, line width=.5pt},
 figblock/.style={draw=black, fill=black!14, line width=.45pt},
 figguide/.style={draw=black!45, line width=.35pt, dash pattern=on 2pt off 2pt},
 figleader/.style={draw=black!60, line width=.35pt}
}

\usepackage{tikz}
\usetikzlibrary{patterns,arrows.meta,patterns.meta,calc}

\definecolor{supportblue}{RGB}{45,75,150}
\definecolor{supportgreen}{RGB}{30,125,85}
\definecolor{supportred}{RGB}{175,50,60}
\tikzset{axe/.style={thin,-{Stealth[length=1.6mm]}},
 original/.style={draw=supportblue,pattern=north east lines,pattern color=supportblue,line width=.4pt},
 enlarged/.style={draw=supportgreen,pattern=north west lines,pattern color=supportgreen,line width=.4pt},
 shifted/.style={draw=supportred,pattern=north east lines,pattern color=supportred,line width=.4pt}}

\setlist[enumerate]{label=\textup{(\arabic*)}}

\numberwithin{equation}{section}
\newtheorem{theorem}{Theorem}[section]
\newtheorem{proposition}[theorem]{Proposition}
\newtheorem{lemma}[theorem]{Lemma}

\newtheorem{corollary}[theorem]{Corollary}

\theoremstyle{definition}

\theoremstyle{remark}
\newtheorem{remark}[theorem]{Remark}

\newcommand{\dst}{\displaystyle}
\newcommand{\cc}{\mathcal C}
\newcommand{\qq}{\mathcal Q}
\newcommand{\R}{\mathbb R}
\newcommand{\Z}{\mathbb Z}
\newcommand{\N}{\mathbb N}
\newcommand{\C}{\mathbb C}
\newcommand{\Q}{\mathbb Q}
\newcommand{\1}{\mathbbm 1}
\newcommand{\one}{\mathbf 1_d}
\newcommand{\Sclass}{\mathcal S}
\newcommand{\Zopen}{\mathcal Z}
\DeclareMathOperator{\supp}{supp}
\DeclareMathOperator{\dist}{dist}
\DeclareMathOperator{\interior}{int}
\DeclareMathOperator{\vol}{vol}
\DeclareMathOperator{\Sp}{Sp}
\DeclareMathOperator{\sinc}{sinc}
\DeclareMathOperator{\sep}{sep_\infty}

\title[Prescribed holes in the time frequency plane]
{Prescribed holes in the time frequency plane and ptychographic phase retrieval}
\author[J. de Dios Pont]{Jaume de Dios Pont}
\address{Center for Data Science, New York University, New York, New York 10011, USA}
\email{jdedios@nyu.edu}

\author[P. Jaming]{Philippe Jaming}
\address{Univ. Bordeaux, CNRS, Bordeaux INP, IMB, UMR 5251, F-33400 Talence, France}
\email{philippe.jaming@math.u-bordeaux.fr}

\author[L. Liehr]{Lukas Liehr}
\address{Department of Mathematics, Bar-Ilan University, Ramat-Gan 5290002, Israel}
\email{lukas.liehr@biu.ac.il}

\author[M. A. Taylor]{Mitchell A. Taylor}
\address{Department of Mathematics, ETH Z\"urich, R\"amistrasse 101,
8092 Z\"urich, Switzerland}
\email{mitchell.taylor@math.ethz.ch}

\date{\today}


\begin{document}

\begin{abstract}
    We investigate the geometry of zero sets of ambiguity functions and its implications for short-time Fourier transform phase retrieval. We prove that every nonzero window with a bounded open zero set does phase retrieval. Moreover, we show that any bounded regular open set, after a sufficiently small rescaling and translation away from the origin, can be realized together with its reflection as the exact open zero set of an ambiguity function. This shows that phase retrieval is possible even when the ambiguity function vanishes on prescribed open sets.
    
    We also establish global geometric obstructions. For example, certain annuli around the origin cannot occur even as subsets of zero sets of Ambiguity functions. These obstructions rely on symplectic Fourier invariance and explicit Laguerre-Gaussian sign certificates. Our three main results are also formally verified in Lean 4.
\end{abstract}

\maketitle
\section{Introduction and results}

\subsection{Problem overview}

In this paper, we investigate two closely related questions in time-frequency
analysis, namely geometric properties of the support of the ambiguity function and
uniqueness in phase retrieval from short-time Fourier transform measurements.
Recall that the short-time Fourier transform (STFT) of a function $f\in L^2(\R^d)$ with
window function $w\in L^2(\R^d)$ is defined by
$$
 V_wf(x,\omega)
 =\int_{\R^d}f(t)\overline{w(t-x)}
   e^{-2\pi i\omega\cdot t}\,dt,
 \qquad (x,\omega)\in\R^d\times\R^d.
$$
The phase retrieval problem for the STFT asks whether $|V_wf|$, the so-called spectrogram, determines
$f$ up to multiplication by a constant of modulus one. Accordingly, we say
that $w$ \emph{does phase retrieval} if, for all $f,g\in L^2(\R^d)$,
the condition
$$
|V_wf|=|V_wg|
$$
implies that $g=e^{i\theta}f$ for some $\theta\in\R$.

\begin{figure}[t]
\vspace*{10pt}
\centering
\begin{tikzpicture}[paperfigure,x=1cm,y=1cm,
  ptyobject/.pic={
    \begin{scope}
      \clip (0,-1.1) rectangle (4.8,1.1);
      \fill[black!30,even odd rule]
        plot[smooth cycle,tension=.75] coordinates {
          (.30,.10) (.45,.65) (.95,.92) (1.45,.68) (1.95,.86) (2.10,.38) (1.72,.02)
          (1.95,-.45) (1.62,-.90) (1.08,-.68) (.62,-.95) (.30,-.50)}
        (.85,.28) ellipse[x radius=.26,y radius=.19]
        plot[smooth cycle,tension=.8] coordinates {
          (2.42,.55) (2.75,.93) (3.25,.80) (3.45,.36) (3.05,.22) (2.85,-.14) (2.50,.06)}
        plot[smooth cycle,tension=.8] coordinates {
          (3.55,-.22) (3.95,.14) (4.50,.28) (4.62,-.30) (4.35,-.85) (3.80,-.92) (3.45,-.60)}
        [rotate around={-18:(2.58,-.70)}] (2.58,-.70) ellipse[x radius=.30,y radius=.17];
      \fill[black!30] (4.28,.76) circle[radius=.15];
    \end{scope}
  }]

  \pic at (0,0) {ptyobject};
  \foreach \c in {1.3,2.4,3.5}{
    \fill[white,fill opacity=.30] (\c,0) circle[radius=.82];
    \fill[black,fill opacity=.07] (\c,0) circle[radius=.82];
  }
  \foreach \c/\cc in {1.3/2.4,2.4/3.5}{
    \begin{scope}
      \clip (\c,0) circle[radius=.82];
      \fill[pattern={Lines[angle=45,distance=2.4pt,line width=.3pt]},
            pattern color=black!60] (\cc,0) circle[radius=.82];
    \end{scope}
  }
  \draw[figoutline] (0,-1.1) rectangle (4.8,1.1);
  \foreach \c in {1.3,3.5}{
    \draw[draw=black!75,line width=.4pt,dash pattern=on 2.2pt off 1.6pt]
      (\c,0) circle[radius=.82];
  }
  \draw[draw=black,line width=.8pt] (2.4,0) circle[radius=.82];
  \foreach \c/\j in {1.3/1,2.4/2,3.5/3}{
    \fill (\c,0) circle[radius=.04];
    \node[below=2pt,inner sep=.8pt,font=\footnotesize] at (\c,0) {$x_\j$};
  }
  \fill[black!10] plot[domain=1.42:3.38,samples=60,smooth]
    (\x,{1.4+.62*exp(-((\x-2.4)/.36)^2)}) -- cycle;
  \draw[draw=black!55,line width=.35pt] (1.42,1.4)--(3.38,1.4);
  \draw[draw=black,line width=.6pt] plot[domain=1.42:3.38,samples=60,smooth]
    (\x,{1.4+.62*exp(-((\x-2.4)/.36)^2)});
  \draw[figguide] (1.58,1.4)--(1.58,1.1);
  \draw[figguide] (3.22,1.4)--(3.22,1.1);
  \node[anchor=west,font=\footnotesize] at (2.95,1.85) {probe $g(\cdot-x_2)$};
  \fill[black!70] (1.85,.34) circle[radius=.025];
  \draw[figleader] (1.85,.34)--(.98,1.42);
  \node[anchor=south,inner sep=1.5pt,font=\footnotesize] at (.72,1.4) {overlap};

  \draw[figaxis,draw=black,line width=.6pt] (4.97,0)--(6.13,0)
    node[midway,above=2pt,font=\footnotesize] {far field};
  \foreach \c/\j/\rot/\ys in {6.85/1/15/.80, 8.2/2/-35/.92, 9.55/3/65/.70}{
    \begin{scope}[shift={(\c,0)}]
      \begin{scope}
        \clip (-.6,-.6) rectangle (.6,.6);
        \begin{scope}[rotate=\rot,yscale=\ys]
          \foreach \rr/\g/\lw in {.17/85/.6,.32/65/.55,.47/48/.5,.62/34/.45,.77/22/.4}{
            \draw[draw=black!\g,line width=\lw pt] (0,0) circle[radius=\rr];
          }
          \foreach \ang/\rad in {25/.245,100/.395,165/.245,235/.395,290/.245,335/.545,60/.545,200/.695}{
            \fill[black!60] (\ang:\rad) circle[radius=.028];
          }
        \end{scope}
        \fill[black!28] (0,0) circle[radius=.1];
        \fill[black!85] (0,0) circle[radius=.052];
      \end{scope}
      \draw[figoutline] (-.6,-.6) rectangle (.6,.6);
      \node[below=2pt,font=\footnotesize] at (0,-.6) {$I(x_\j,\cdot)$};
    \end{scope}
  }
  \draw[draw=black,line width=.8pt] (7.6,-.6) rectangle (8.8,.6);

  \draw[figaxis,draw=black,line width=.6pt] (10.32,0)--(12.28,0)
    node[midway,above=2pt,font=\footnotesize] {phase retrieval};
  \pic[scale=.58] at (12.45,0) {ptyobject};
  \draw[figoutline] (12.45,-.638) rectangle (15.234,.638);

  \node[font=\footnotesize] at (2.4,-1.45) {object $f$ and probe positions};
  \node[font=\footnotesize] at (8.2,-1.45) {diffraction patterns};
  \node[font=\footnotesize] at (13.84,-1.45) {reconstruction of $f$};
\end{tikzpicture}
\caption{Schematic illustration of ptychographic imaging. A localized probe $g$ is moved to overlapping positions $x_1,x_2,x_3$ on the object $f$. At each position $x_j$, the detector records the far-field intensity $I(x_j,\cdot)=|V_gf(x_j,\cdot)|^2$. The redundancy provided by the overlap allows phase retrieval, that is, the reconstruction of $f$ up to a constant phase factor.}
\label{fig:ptychography}
\end{figure}
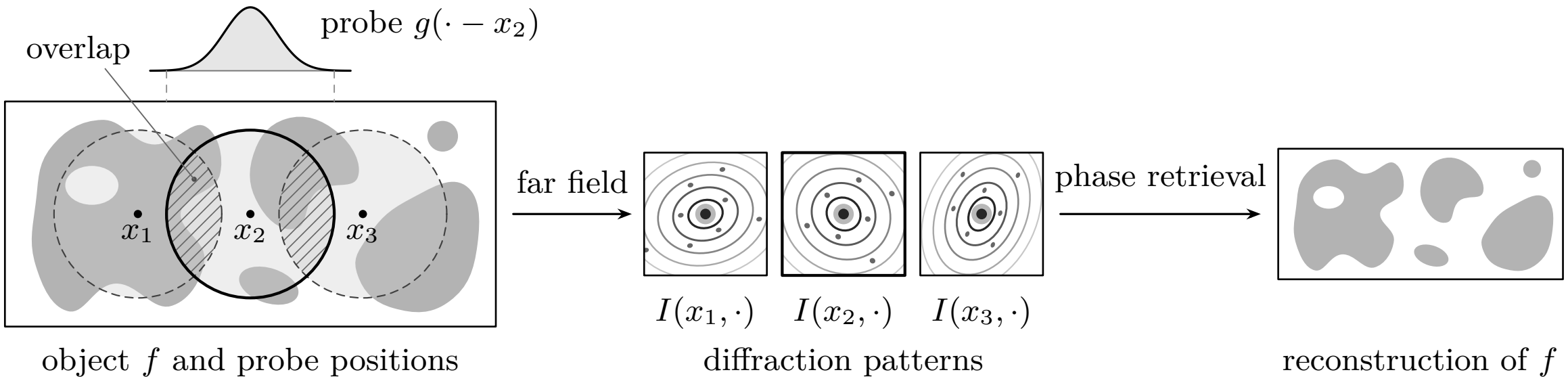

Questions of this type arise in numerous physical applications, most notably
ptychography \cite{Marchesini2016,Rodenburg2008,Shechtman2015}.
Ptychography is a coherent imaging technique that aims to recover
both the amplitude and the phase of an unknown object from
intensity-only diffraction measurements. Originally developed
in electron microscopy, it is now widely used in X-ray imaging,
optical microscopy, and materials science.
In a typical ptychographic experiment, a localized illumination
beam, called the \emph{probe}, is scanned across the object,
and the intensity of the resulting far-field diffraction pattern
is recorded at each position (see Figure~\ref{fig:ptychography}). More precisely, if
$f:\mathbb R^d\to\mathbb C$ represents the object and
$g:\mathbb R^d\to\mathbb C$ the probe, the measurements take
the form
\[
I(x,\omega)
=
\left|
\int_{\mathbb R^d}
f(t)\overline{g(t-x)}
e^{-2\pi i\omega\cdot t}\,\mathrm{d}t
\right|^2
=
|V_gf(x,\omega)|^2.
\]
In practice, measurements are collected at a discrete set of scanning positions and detector frequencies. The fundamental principle of ptychography is that consecutive illumination regions overlap, introducing redundancy into the
measurements that can be exploited to recover the missing phase. Thus, from a mathematical perspective, ptychography can be viewed
as a structured STFT phase retrieval problem, in which one seeks
to reconstruct $f$, up to a global phase factor, from partial
measurements of $|V_gf|$.
The practical importance of ptychography makes it an important application
of time-frequency analysis and  has therefore
attracted considerable attention in recent years \cite{AlaifariWellershoff2021,Bartusel2023,GrohsKoppensteinerRathmair2020}.

\subsection{The ambiguity function}

A standard sufficient condition for a function $w$ to do phase retrieval is
expressed in terms of its ambiguity function. For $f,g\in L^2(\R^d)$, we
define
\[
 A(f,g)(x,\omega)
 :=\int_{\R^d}
 f\!\left(t+\frac{x}{2}\right)
 \overline{g\!\left(t-\frac{x}{2}\right)}
 e^{-2\pi i\omega\cdot t}\,dt,
 \qquad (x,\omega)\in\R^d\times\R^d,
\]
and call $A(f):=A(f,f)$ the ambiguity function of $f$. 
The ambiguity function and the short-time Fourier transform are basic
objects in time-frequency analysis. Both are used to represent the function $f$ in
the time-frequency plane, with the delay-variable $x$ representing time and $\omega$
representing the frequency component of $f$.
We refer to
\cite{Folland1989,Grochenig2001} and Section~\ref{sec:preliminaries}
for their properties.

The spectrogram $|V_wf|^2$ and the ambiguity function are linked by the
identity
\[
 \mathcal F_\sigma(|V_wf|^2)=A(f)\overline{A(w)},
\]
where $\mathcal F_\sigma$ denotes the symplectic Fourier transform (see
Section~\ref{sec:preliminaries}). This identity makes it clear why the zeros of
$A(w)$ play a central role in ptychographic phase retrieval. If $A(w)$ never vanishes, equality of two STFT
magnitudes gives $A(f)=A(g)$, which in turn implies that
$f=e^{i\theta}g$ for some $\theta\in\R$. Thus, every window whose ambiguity
function is nowhere zero does phase retrieval. The main question is how much information can be lost while
maintaining the phase retrieval property. More precisely, for which $\Omega=\overline{\text{supp}A(w)}$ is it true that $A(f)=A(g)$ on $\Omega$ still implies that 
$f=e^{i\theta}g$ for some $\theta\in\R$?

Gr\"ochenig, Jaming, and Malinnikova
\cite{GrochenigJamingMalinnikova2020} studied the problem of determining
when the zero set of a Wigner distribution, short-time Fourier transform, or
ambiguity function is \emph{empty}. In particular, they identified non-Gaussian functions, including one-sided exponentials and certain convolutions, whose ambiguity functions are zero-free. Recall also
Hudson's theorem \cite{Hudson1974}, which characterizes the functions with
nonnegative Wigner distribution as the generalized Gaussians.

Being zero-free is, however, not necessary for phase retrieval. As the ambiguity function is continuous, phase retrieval is also implied whenever $A(w)$ is nonzero on a dense set which is by now standard result in the phase retrieval literature
\cite{AlaifariWellershoff2021,BojarovskaFlinth2016,GrohsKoppensteinerRathmair2020,GrohsRathmair2019}.

\smallskip

Windows that violate the property of being nonzero on a dense set arise naturally. To see this, suppose that
$w$ is supported in a compact set $K\subset\R^d$. If $x\notin K-K$, then
$w$ and its translate $w(\cdot-x)$ have disjoint supports, and the definition of the ambiguity function implies that $A(w)(x,\cdot)$ vanishes identically. Hence
\[
\supp A(w)\subset (K-K)\times\R^d,
\]
and the zero set of $A(w)$ contains the unbounded open set
$(\R^d\setminus(K-K))\times\R^d$.

Windows of this type {\em do not} do phase retrieval. Indeed, choose
$f_1,f_2\in L^2(\R^d)\setminus\{0\}$ with bounded and sufficiently separated
supports so that the support of no translate of $w$ intersects both support of $f_1$ and $f_2$ simultaneously. Then, at each point
of $\R^{2d}$, at most one of $V_wf_1$ and $V_wf_2$ is nonzero which implies
that $|V_w(f_1+e^{i\theta}f_2)|$
does not depend on $\theta\in\R$, whereas $f_1+e^{i\theta}f_2$ and
$f_1+f_2$ do not agree up to a global phase when
$\theta\notin2\pi\Z$. Ambiguities of this kind for functions with disconnected
support are well known
\cite{Bartusel2023,GrohsRathmair2019}.
In all these examples of windows which do not do phase retrieval, the open
set on which the ambiguity function vanishes is unbounded.
It is therefore natural to ask how much vanishing of $A(w)$ can be tolerated without losing injectivity. This question has been raised explicitly in the recent literature. Bartusel formulated a closely related problem in \cite[Remark~5.4.1]{Bartusel2023}, and later constructed, in every dimension, a window $w$ such that $\{A(w)\neq0\}$ is not of full measure but nevertheless dense, thus still doing phase retrieval \cite[Example~8.2.2]{Bartusel2024}.  It is therefore useful to introduce the open zero set
\[
 \Zopen(w):=\R^{2d}\setminus\supp A(w)
 =\interior\{z\in\R^{2d}:A(w)(z)=0\}.
\]
We note that the windows in Bartusel's example satisfy $\Zopen(w)=\varnothing$. See also \cite{Bartusel2026} for quantitative
refinements on discrete groups, \cite{AlaifariYang2025} for a multiwindow
approach designed to compensate for zeros of individual ambiguity functions
and the problem collection \cite{Phasebook2025}.

\subsection{Main results}

Our first result shows that the ambiguity function may vanish on a bounded open set while the corresponding window still does phase retrieval. More precisely, the boundedness of the open zero set $\Zopen(w)$ is sufficient for $w$ to do phase retrieval.

\begin{theorem}
\label{thm:phase-retrieval-d}
Let $w\in L^2(\R^d) \setminus \{ 0 \}$.  If $\Zopen(w)$ is bounded then $w$ does phase retrieval.
\end{theorem}

Theorem~\ref{thm:phase-retrieval-d} shows that the full-support criterion is
far from necessary: as long as the open zero set $\Zopen(w)$ is bounded, the
ambiguity function of $w$ may vanish on an open set without affecting phase
retrieval. In other words, phase retrieval can only fail if $\Zopen(w)$ is
unbounded, as in the examples above.
In Remark~\ref{rem:finite-measure}, we show that the boundedness of $\Zopen(w)$ can even be replaced by the assumption that $\Zopen(w)$ has finite measure.

\smallskip

Theorem~\ref{thm:phase-retrieval-d} turns a question about injectivity into
a geometric one:

\medskip

\centerline{\it Which open subsets of phase space can occur as $\Zopen(w)$?}

\medskip

\noindent This is the central question of the present paper.
As mentioned above, it is easy to construct functions $w$ for which $\Zopen(w)$ is nonempty.
For instance, if the support $S$ of $w$ is such that $\R^d\setminus(S-S)$ contains an open interval $I$, then
\[
I \times \R^d \subset \Zopen(w)
\]
so that $\Zopen(w)$ contains a strip. Further examples may be obtained by
applying a fractional Fourier transform to $w$, whose effect is to rotate the support \cite{Wilcox1991}. Further, Melenk and Zimmermann \cite{MelenkZimmermann} constructed functions $w$ such that both the support $S$ of $w$ and the support $\widehat{S}$ of its Fourier transform are such that $\R^d\setminus(S-S)$ and  $\R^d\setminus(\widehat{S}-\widehat{S})$ contain open intervals $I$ and $J$, respectively, leading to the inclusion
\[
(I \times \R^d) \cup (\R^d \times J) \subset \Zopen(w).
\]
All these examples, however,
produce unbounded open zero sets and hence do not fall within the scope of
Theorem~\ref{thm:phase-retrieval-d}.

This leads us to investigate more systematically the possible open zero sets
of ambiguity functions. There are first some elementary restrictions. Since
$A(w)$ is continuous on $\R^{2d}$ and satisfies
\[
A(w)(0,0)=\|w\|_{L^2(\R^d)}^2\neq0
\quad \text{and} \quad
A(w)(x,\omega)=\overline{A(w)(-x,-\omega)},
\]
the set $\Zopen(w)$ is regular open, does not contain the origin, and
is centrally symmetric:
\[
0\notin\Zopen(w),
\qquad
\Zopen(w)=-\Zopen(w).
\]
These restrictions follow directly from the continuity and the symmetry of the
ambiguity function. Further, less obvious restrictions come from uncertainty
principles, which assert that a nonzero function and its Fourier transform
cannot both be too concentrated. Since the ambiguity function encodes the time
and the frequency content of a function simultaneously, such principles limit
the size of the open zero set of an ambiguity function. Among the various mathematical formulations of the
uncertainty principle, the one
that is of particular interest here is the classical Amrein--Berthier \cite{AmreinBerthier} and
Benedicks \cite{Benedicks1985} theorem for the Fourier transform:

\vspace{0.2cm}
\begin{center}
    \textit{A non-zero function $f\in L^2(\R^d)$ and its Fourier transform
$\widehat{f}$ can not both have \\ a support of finite measure.}
\end{center}
\vspace{0.2cm}

Observe that an analogous result holds for the ambiguity function
\cite{JamingCRAS,Janssen1998,Wilczok}: if
$f\in L^2(\R^d)$ is nonzero, then $\supp A(f)$ has infinite measure (see also \cite{BonamiDemange,BonamiDemangeJaming2003,Demange2005} for further statements in this direction). In particular,
$\Zopen(f)$ cannot be the complement of a set of finite measure.

A related result of Amrein and Berthier \cite{AmreinBerthier} states that for two sets $S,\Sigma\subset\R^d$ of finite
measure, the space
$$
\mathcal{F}(S,\Sigma)
=
\left\{
f\in L^2(\R^d)\,:\ 
\supp f\subset\R^d\setminus S,\ 
\supp\widehat{f}\subset\R^d\setminus\Sigma
\right\}
$$
is infinite dimensional. This raises the question of whether a comparable flexibility holds for ambiguity functions. We show that, locally, the answer is affirmative: sufficiently small open sets satisfying the elementary geometric restrictions above can be prescribed as holes in the support of an ambiguity function. Globally, however, additional obstructions arise. Nevertheless, the following holds:

\begin{theorem}
\label{thm:main}
Let $\Omega\subset\R^{2d}$ be bounded, nonempty, and regular
open, and let $a\in\R^{2d}\setminus\{0\}$.  There exists
$b_0=b_0(d,\Omega,a)>0$ such that for every $0<b<b_0$ one can find
$f\in L^2(\R^d)$ with
\[
 \Zopen(f)
 =(a+b\Omega)\cup (-(a+b\Omega) ).
\]
\end{theorem}

\begin{figure}[htbp]
\centering
\begin{tikzpicture}[paperfigure,x=1.3cm,y=1.3cm,
  omega/.pic={
    \path[pic actions,even odd rule]
      plot[smooth cycle,tension=.7] coordinates {
        (-1.55,-.05) (-1.35,.55) (-.85,.95) (-.35,.70) (.05,.95) (.55,1.00)
        (.85,.60) (.60,.20) (.95,-.15) (.80,-.70) (.30,-1.00) (-.15,-.70)
        (-.60,-.95) (-1.15,-.80) (-1.30,-.45)}
      (-.78,.12) ellipse[x radius=.27,y radius=.20];
    \path[pic actions] (1.45,.66) ellipse[x radius=.22,y radius=.17];
    \path[pic actions] plot[smooth cycle,tension=.8] coordinates {
      (1.20,-.60) (1.36,-.36) (1.62,-.44) (1.56,-.76) (1.34,-.82)};
    \path[pic actions] (-1.78,.78) circle[radius=.10];
  }]
  \def\figax{3.7}\def\figay{.9}\def\figb{.35}
  \pic[fill=black!12,draw=black!45,line width=.7pt] at (0,0) {omega};
  \draw[figaxis] (-5.1,0)--(5.15,0) node[right] {$x$};
  \draw[figaxis] (0,-1.75)--(0,1.9) node[above] {$\omega$};
  \fill (0,0) circle[radius=.028];
  \node[below left,inner sep=2pt,font=\footnotesize] at (0,0) {$0$};
  \node at (.22,.52) {$\Omega$};
  \pic[scale=\figb,fill=supportblue!28,draw=supportblue,line width=.6pt]
    at (\figax,\figay) {omega};
  \pic[scale=\figb,rotate=180,fill=supportblue!28,draw=supportblue,line width=.6pt]
    at (-\figax,-\figay) {omega};
  \node[above=2pt,text=supportblue!85!black] at (\figax,\figay+.40) {$a+b\Omega$};
  \node[below=2pt,text=supportblue!85!black] at (-\figax,-\figay-.40) {$-(a+b\Omega)$};
  \draw[figaxis,draw=black!75] (.95,1.08) to[out=40,in=160]
    node[pos=.5,above=1pt,font=\footnotesize] {$z\mapsto a+bz$} (\figax-.72,\figay+.22);
  \draw[figaxis,draw=black!75] (-.95,-1.08) to[out=220,in=340]
    node[pos=.5,below=1pt,font=\footnotesize] {$z\mapsto -(a+bz)$} (-\figax+.72,-\figay-.22);
\end{tikzpicture}
\caption{Illustration of Theorem~\ref{thm:main} in the time-frequency plane $\R^2$. A bounded regular open set $\Omega$ (grey) is scaled by a factor $b>0$ and translated by $a\neq0$. If $b<b_0$, then the resulting set together with its reflection through the origin is the open zero set $\Zopen(f)=(a+b\Omega)\cup(-(a+b\Omega))$ of an ambiguity function (blue).}
\label{fig:main-theorem}
\end{figure}
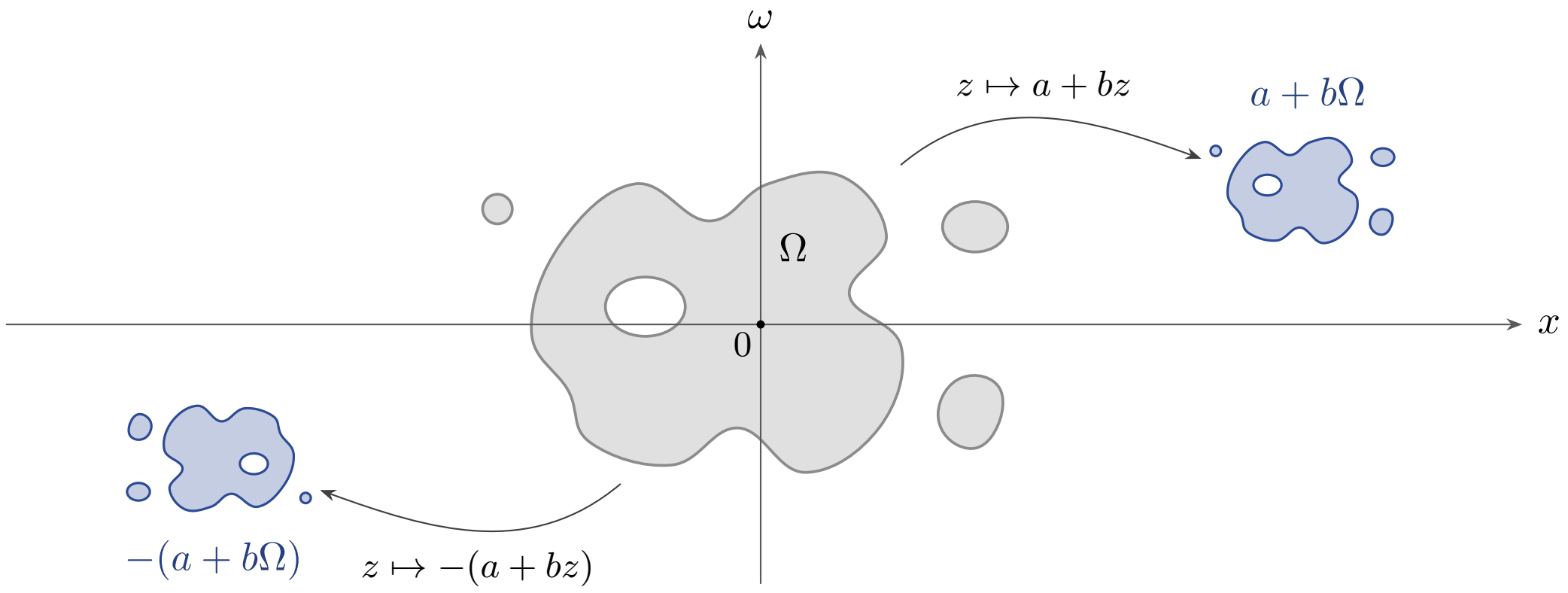

Thus any bounded regular open set can occur, after rescaling,
together with the centrally reflected copy
as the open zero set of an ambiguity function (see Figure~\ref{fig:main-theorem}). Since the resulting open
zero set is bounded, it follows from Theorem~\ref{thm:phase-retrieval-d} that
every window supplied by Theorem~\ref{thm:main} also does phase retrieval.
In particular, phase retrieval is
compatible with genuine open gaps in the support of the ambiguity function.
We also remark that a quantitative version of Theorem~\ref{thm:main}, with an explicit sufficient bound on the scaling parameter $b_0$, is given in Corollary~\ref{cor:quantitative}.

Globally, however, the elementary restrictions above do not suffice. To see
this, define for $0<r<R<\infty$ the annulus
\begin{equation}
 O_{r,R}:=\{z\in\R^{2d}:r<|z|<R\}.
 \label{eq:annulus}
\end{equation}
The annulus $O_{r,R}$ is bounded, regular open, centrally symmetric, and does
not contain the origin, so it satisfies all necessary restrictions to be the open zero set of an ambiguity function.
Nevertheless, the next theorem shows that, for suitable radii, $O_{r,R}$
cannot even be contained in the open zero set of a nonzero ambiguity function.

\begin{theorem}
\label{thm:counterexample}
There exist $0 < r < R < \infty$ such that no nonzero $f\in L^2(\R)$ satisfies
$$
 O_{r,R}\subset\Zopen(f).
$$
\end{theorem}

Symplectic images of the annulus $O_{r,R}$ from Theorem~\ref{thm:counterexample} and higher-dimensional analogues are treated in Section~\ref{sec:annuli}, see Corollaries~\ref{cor:symplectic} and~\ref{cor:all-dimensional-annuli}.

\subsection{Outline of the proofs}

We briefly describe the main ideas behind the proofs.

To prove Theorem \ref{thm:main}, we first place the center of the prescribed
hole on the delay axis, which in the following refers to the set $\{(x,0):x\in\R^d\}\subset\R^{2d}$, and construct a reference function as a sparse sum of characteristic functions of cubes. Support considerations are arranged so that, for delays near the two selected points, the relevant part of the ambiguity function depends linearly on the perturbation.

The frequency variable is controlled through a non-harmonic Fourier series whose frequencies lie in the union of a lattice and an irrational translate of that lattice. The presence of the second translate is essential, as a purely periodic Fourier series would reproduce its zero set in every period, whereas the irrational translate allows us to prescribe a single bounded zero set. Suitable linear combinations of translates of ambiguity functions of characteristic functions of cubes (see Lemma~\ref{lem:one-coordinate-difference}) then provide the
required localization in the delay variable $x$.

A countable perturbation argument makes the ambiguity function
nonzero on a dense subset of the complement of the prescribed set in
a strip $\Delta\times\R^d$, where $\Delta\subset\R^d$ is the union of two small cubes. A second perturbation argument fills the complement outside this strip while leaving the already constructed part unchanged. Finally, symplectic covariance moves the construction from the delay axis to an arbitrary prescribed nonzero center.

The obstruction in Theorem~\ref{thm:counterexample} has a quite different
origin. Since for every $f\in L^2(\R^d)$ we have
$\mathcal F_\sigma(|A(f)|^2)=|A(f)|^2$,
it follows that $|A(f)|^2$ is orthogonal to every real Schwartz function $h$
satisfying $\mathcal F_\sigma h=-h$.
If $h(0)>0$ and the set $\{z:h(z)<0\}$ were contained in $\Zopen(f)$,
then $|A(f)|^2h$ would be nonnegative everywhere and strictly positive near
the origin, contradicting this orthogonality. An explicit combination of
Laguerre--Gaussian functions produces such a function $h$ for which the set $\{z:h(z)<0\}$
is precisely an annulus $O_{r,R}$ as in Theorem~\ref{thm:counterexample}. The same method also yields
higher-dimensional examples.

Sign certificates of this kind are familiar in Fourier optimization, most
prominently in the linear programming bounds for sphere packings
\cite{CohnElkies2003} and in the sign uncertainty principles of Bourgain,
Clozel, and Kahane \cite{BourgainClozelKahane2010}. We refer to
Remark~\ref{rem:sign-uncertainty} for further details.

Finally, the proof of Theorem~\ref{thm:phase-retrieval-d} relies on a
rigidity property of ambiguity functions. We show that if the difference of two ambiguity functions has compact support, then the underlying functions must be linearly dependent.
Using that $\Zopen(w)$ is bounded, we can apply Lemma~\ref{lem:compact-rigidity} and obtain phase retrieval.

\section{Preliminaries on the ambiguity function}
\label{sec:preliminaries}

\subsection{Basic facts}
We start by listing basic facts on the ambiguity function that can be found for instance
in \cite{AuslanderTolimieri1985,Folland1989,Grochenig2001,Wilcox1991}. We adopt here the choice of normalization
\[
 A(f,g)(x,\omega)
 :=\int_{\R^d}
 f\!\left(t+\frac{x}{2}\right)
 \overline{g\!\left(t-\frac{x}{2}\right)}
 e^{-2\pi i\omega\cdot t}\,dt,
 \qquad (x,\omega)\in\R^d\times\R^d,
\]
and let $A(f)=A(f,f)$ be the ambiguity function of $f$.
It will often be convenient to write $z=x+i\omega$ and $A(f,g)(z)=A(f,g)(x,\omega)$.
The variable $x$ will be called the {\em delay-variable}.

Recall that, for $f\in L^p(\mathbb{R}^d)$, we define its essential support by
\[
\operatorname{supp}f
=\mathbb{R}^d\setminus
\bigcup\left\{U\subset\mathbb R^d:
U\text{ open},\ f=0\text{ a.e. on }U\right\}
\]
which is a closed set.
For sets $B,B'\subset\R^d$ and $s\in\R^d$, we write $B+s=\{t+s:t\in B\}$ and
$B-B'=\{t-t':t\in B,\ t'\in B'\}$. 
Recall that if $B$ and $B'$ are closed
(but none being compact), then
$B-B'$ needs not be closed. Moreover, $\tau_sf(t)=f(t-s)$ denotes the
translation of $f$ by $s$, and $u_+=\max(u,0)$ for $u\in\R$. The sinc function
is defined by $\sinc x=\frac{\sin x}{x}$ for $x\neq0$ and $\sinc0=1$. Finally, $\mathcal C_0(\R^{2d})$ denotes the space of continuous functions on $\R^{2d}$ that vanish at infinity.

\begin{lemma}\label{lem:basic-facts}
Let $f,g\in L^2(\R^d)$. Then the following statements hold.
\begin{enumerate}
\item\label{fact:continuity}
$A(f,g)\in\mathcal C_0(\R^{2d})$ and $\|A(f,g)\|_\infty\leq\|f\|_2\|g\|_2$.
Moreover, the map $(f,g)\mapsto A(f,g)$ from $L^2(\R^d)\times L^2(\R^d)$ to
$\mathcal C_0(\R^{2d})$ is continuous, linear in $f$ and anti-linear in $g$.
\item\label{fact:uniqueness}
$A(f)=A(g)$ if and only if there exists $\theta\in\R$ such that
$f=e^{i\theta}g$.
\item\label{fact:sym}
$A(g,f)(x,\omega)=\overline{A(f,g)(-x,-\omega)}$ for all
$(x,\omega)\in\R^{2d}$. In particular, $A(f)(-z)=\overline{A(f)(z)}$ for all
$z\in\R^{2d}$, and $\|A(f)\|_\infty=A(f)(0)=\|f\|_2^2$.
\item\label{fact:obstruction}
If $f\neq0$, then the zero set
\[
Z_f=\{z\in\R^{2d}:A(f)(z)=0\}
\]
is closed and satisfies $Z_f=-Z_f$ and $0\notin Z_f$. Consequently, the open
zero set $\Zopen(f)=\interior Z_f$ is regular open and satisfies
$\Zopen(f)=-\Zopen(f)$.
\item\label{fact:support}
$\supp A(f,g)\subset\overline{\supp f-\supp g}\times\R^d$. In particular,
\begin{equation}
 \supp_x A(f,g)\subset \overline{\supp f-\supp g},
 \label{eq:support-rule}
\end{equation}
where $\supp_x$ denotes the support in the delay variable, that is,
$A(f,g)(x,\cdot)$ vanishes identically whenever $x\notin\supp f-\supp g$.
If $f$ and $g$ are nonzero and compactly supported, then
$\supp A(f,g)=(\supp f-\supp g)\times\R^d$.
\item\label{fact:translation}
For all $r,s\in\R^d$ and $(x,\omega)\in\R^{2d}$, one has
\[
A(\tau_rf,\tau_sg)(x,\omega)=e^{-\pi i\omega\cdot(r+s)}A(f,g)(x-(r-s),\omega).
\]
\item\label{fact:expra1q}
For $L>0$ and $\qq=[0,1]^d$,
\[
A(\1_{L\qq})(x,\omega)=\prod_{j=1}^d e^{-\pi iL\omega_j}(L-|x_j|)_+\sinc\bigl(\pi(L-|x_j|)\omega_j\bigr).
\]
\end{enumerate}
\end{lemma}

Property~\ref{fact:uniqueness} is crucial for phase retrieval. Property~\ref{fact:continuity} implies that, for every $z\in\R^{2d}$, the map
\[
s\mapsto A(f+sg)(z)=A(f)(z)+s\bigl(A(f,g)(z)+A(g,f)(z)\bigr)+s^2A(g)(z)
\]
is a polynomial of degree at most $2$ on $\R$.

\subsection{The symplectic Fourier transform}

Write the standard symplectic form on $\R^{2d}$ as
\[
 \sigma((x,\omega),(y,\eta))
 =\omega\cdot y-x\cdot\eta=(x,\omega)\cdot(-\eta,y).
\]
We use the symplectic Fourier transform
\[
 (\mathcal F_\sigma F)(x,\omega)
 =\iint_{\R^{2d}}F(y,\eta)
   e^{-2\pi i\sigma((x,\omega),(y,\eta))}\,dy\,d\eta,
 \qquad F\in L^1(\R^{2d}).
\]
It arises from the Euclidean Fourier transform on $\R^{2d}$ by the linear
change of variables $(y,\eta)\mapsto(\eta,-y)$. In particular,
$\mathcal F_\sigma$ extends to a unitary involution of $L^2(\R^{2d})$.
Recall also that the orthogonality relations of the short-time
Fourier transform \cite[Chapter~3]{Grochenig2001} imply that
$\|V_wf\|_{L^2(\R^{2d})}=\|f\|_2\|w\|_2$ for all
$f,w\in L^2(\R^d)$. In particular, we have
$|V_wf|^2\in L^1(\R^{2d})$.

The spectrogram and the ambiguity function
are related by the following identity.

\begin{lemma}
\label{lem:spectrogram-identity}
For $f,w\in L^2(\R^d)$, we have $\mathcal F_\sigma(|V_wf|^2)=A(f)\,\overline{A(w)}$.
\end{lemma}

The latter lemma describes the classical ambiguity function relation of
time--frequency analysis. It is recorded, for instance, in
\cite[eq.~(4)]{Bartusel2023}, a version for tempered distributions in
dimension one appears in \cite[Theorem~A.2]{GrohsRathmair2019}, and
classical treatments can be found in
\cite{ClaasenMecklenbrauker1980,Cohen1989}.

The formula was also re-discovered in \cite{JamingCRAS,Janssen1998} where it was used to prove
that the support of $A(f,g)$ has infinite measure.

\subsection{Symplectic covariance}
\label{sec:symplectic}

For $B\in\operatorname{GL}(d,\R)$ and a real symmetric matrix $M$, define
\begin{align*}
 (\mathcal D_Bf)(t)&=|\det B|^{1/2}f(Bt),\\
 (\mathcal C_Mf)(t)&=e^{\pi i t\cdot Mt}f(t),\\
 \mathcal Ff&=\widehat f.
\end{align*}
Direct changes of variables give
\begin{align}
 A(\mathcal D_Bf)(x,\omega)
 &=A(f)(Bx,B^{-T}\omega),
 \label{eq:cov-dilation}\\
 A(\mathcal C_Mf)(x,\omega)
 &=A(f)(x,\omega-Mx),
 \label{eq:cov-chirp}\\
 A(\widehat f)(x,\omega)
 &=A(f)(-\omega,x).
 \label{eq:cov-fourier}
\end{align}
The corresponding transformations of phase space, together with their
inverses, generate the symplectic group $\Sp(2d,\R)$ of all matrices
$S\in\R^{2d\times2d}$ satisfying $\sigma(Sz,Sz')=\sigma(z,z')$ for all
$z,z'\in\R^{2d}$ (see \cite[Chapter~4]{Folland1989} and \cite[Chapter~9]{Grochenig2001}).
Consequently, for every
$S\in\Sp(2d,\R)$ there is an operator $\mathcal U_S$, which is a composition of dilations,
chirps, Fourier transforms, and their inverses, such that
\begin{equation}
 A(\mathcal U_Sf)(z)=A(f)(S^{-1}z),
 \qquad
 \Zopen(\mathcal U_Sf)=S\Zopen(f).
 \label{eq:full-symplectic-covariance}
\end{equation}

\subsection{A perturbation lemma}

We next state a perturbation lemma, which is used in the proof of Theorem~\ref{thm:main}.

\begin{lemma}
\label{lem:diagonal}
Let $E\subset\R^{2d}$, let $h,g_1,g_2,\ldots\in L^2(\R^d)$ be non-zero functions and let
$z_1,z_2,\ldots\in\R^{2d}$. Suppose that the following two conditions hold:
\begin{enumerate}
\item\label{hyp:diagonal-1} for every $N$ and every choice of real numbers $s_1,\ldots,s_N$, we have
      \[
       A\left(h+\sum_{j=1}^Ns_jg_j\right)(z)=0
       \qquad\text{for all }z\in E;
      \]
\item\label{hyp:diagonal-2} for every $N$, if $g \in \operatorname{span}\{ g_j : j < N \}$ satisfies $A(h+g)(z_k) \neq 0$ for all $k < N$, then the polynomial
$$
s \longmapsto A(h+g+sg_N)(z_N)
$$
is not identically zero.
\end{enumerate}
Then there are real numbers $s_n$ for which $f=h+\sum_ns_ng_n$ converges in $L^2(\R^d)$
and such that
\[
 A(f)(z)=0\quad\text{for all }z\in E,
 \qquad
 A(f)(z_n)\neq0\quad\text{for every }n\ge1.
\]
\end{lemma}

For $N=0$, the condition $g\in\operatorname{span}\{g_j\,:\ j<N\}$ of course means that $g=0$
and there is no restriction.

\begin{proof}
We construct the numbers $s_n$ by induction and start by setting $f_0=h$. Since $P_1\,:s\to A(f_0+sg_1)(z_1)$ is a non-zero polynomial of degree
at most $2$, we can choose a positive real number $s_1$ that is not a root of $P_1$ and satisfies $s_1\|g_1\|\leq 1/2$,
and set $f_1=h+s_1g_1$. Then $A(f_1)(z_1)\not=0$ and hypothesis~\ref{hyp:diagonal-1} shows that
$A(f_1)(z)=0$ for all $z\in E$.

Suppose now that $s_1,\ldots,s_N$ have been chosen so that $f_N=h+\sum_{j=1}^Ns_jg_j$
satisfies $A(f_N)(z)=0$ for all $z\in E$ and $A(f_N)(z_j)\not=0$ for $j=1,\ldots,N$.
Define, for $n\leq N$, the quantities $m_n$ and $\delta_n$ via
\[
 m_n:=\min_{1\le j\le n}|A(f_n)(z_j)|>0,
 \qquad
 \delta_n:=\min\left\{1,
   \frac{m_n}{4(\|f_n\|_2+1)}\right\}.
\]

For each $j\le N+1$, the map 
\[
A(f_N+sg_{N+1})(z_j)=A(f_N)(z_j)+s\bigl(A(f_N,g_{N+1})+A(g_{N+1},f_N)\bigr)(z_j)+s^2A(g_{N+1})(z_j)
\]
is a polynomial in $s$ of degree at most two. Further, the induction hypothesis asserts
that this polynomial is not $0$ when $j\leq N$, while hypothesis~\ref{hyp:diagonal-2} asserts that this
is also true when $j=N+1$. It follows that each of those polynomials has at most $2$
roots. In particular, there is an $\eta_{N+1}>0$ such that, if $0<s<\eta_{N+1}$,
then $A(f_N+sg_{N+1})(z_j)$ avoids all roots, {\it i.e.} $A(f_N+sg_{N+1})(z_j)\not=0$ for all $j\leq N+1$.

We may thus choose $s_{N+1}$ so that it satisfies the bound
\begin{equation}
 \|s_{N+1}g_{N+1}\|_2\le 2^{-(N+1)}\min\{1,\delta_1,\ldots,\delta_{N}\}
 \label{eq:diagonal-small}
\end{equation}
and set $\dst f_{N+1}=f_N+s_{N+1}g_{N+1}=h+\sum_{n=1}^{N+1}s_ng_n$.

\smallskip

This finishes the inductive construction of the sequence $(s_n)_{n\geq 1}$.
From \eqref{eq:diagonal-small} we immediately get that the series
$\sum s_ng_n$ converges in $L^2$
so that we can set
$$
f=\lim f_n=h+\sum_{n=1}^{+\infty}s_ng_n.
$$

It remains to show that $A(f)$ has the desired properties.

From the continuity of the ambiguity function, it follows that $A(f_n)\to A(f)$ uniformly over $\R^{2d}$. On the other hand, for
$z\in E$, we have $A(f_n)(z)=0$ and therefore $A(f)(z)=0$ as well.

Next, set $r_n=f-f_n$. Using \eqref{eq:diagonal-small}, we have $\|r_n\|_2<\delta_n$.
Moreover,
\begin{equation}
\begin{split}
|A(f)(z)-A(f_n)(z)|
&= |A(f_n+r_n)(z)-A(f_n)(z)| \\
&= |A(f_n,r_n)(z)+A(r_n,f_n)(z)+A(r_n)(z)| \\
&\leq 2\|f_n\|_2\|r_n\|_2+\|r_n\|_2^2
< 2(\|f_n\|_2+1)\delta_n \\
&\leq \frac{m_n}{2},
\end{split}
\end{equation}
which gives 
$$
|A(f)(z_n)|\geq |A(f_n)(z_n)|-\frac{m_n}{2}\geq \frac{m_n}{2}>0
$$
showing that $A(f)(z_n)$ does not vanish on the sequence $z_n$.
\end{proof}

\section{Non-harmonic Fourier series with prescribed zeros}
\label{sec:fourier-series}

The goal of this section is to construct Fourier series whose zero set is a prescribed compact set $F\subset\R^d$ contained in a translate of $(0,1)^d$ (see Proposition~\ref{prop:exact-gap}). Since the zero set of a $\Z^d$-periodic function is invariant under translations by $\Z^d$, we cannot choose all frequencies of such a series in $\Z^d$, and instead choose them in $\Z^d\cup(\Z^d+\vartheta)$ for a suitable vector $\vartheta\in\R^d$.

Fix a vector $\vartheta\in(0,1/2)^d$ such that
\[
 1,\vartheta_1,\ldots,\vartheta_d
 \quad\text{are linearly independent over }\Q
\]
and define
\[
 \Lambda=\Z^d\cup (\Z^d+\vartheta).
\]
For $x\in\R^d$, $E\subset\R^d$, and a discrete set
$\Lambda\subset\R^d$, define the distance of $x$ to $E$ and the separation constant of $\Lambda$ via
\[
 \dist_\infty(x,E):=\inf_{y\in E}\|x-y\|_\infty,
 \qquad
 \operatorname{sep}_\infty(\Lambda)
 :=\inf_{\substack{\lambda,\mu\in\Lambda\\\lambda\neq\mu}}
 \|\lambda-\mu\|_\infty.
\]
The rational independence of $1,\vartheta_1,\ldots,\vartheta_d$  guarantees that the two cosets $\Z^d$ and $\Z^d+\vartheta$ are disjoint, and the distinct points of $\Lambda$ are uniformly separated with
\[
\operatorname{sep}_\infty(\Lambda)=\dist_\infty(\vartheta,\Z^d)<\dfrac{1}{2}.
\]
We first require the following standard result \cite[Theorem 2.29]{lee2012smooth}.

\begin{lemma}
\label{lem:smooth-peak}
If $F\subset(0,1)^d$ is compact, then there exists
$\eta\in C_c^\infty((0,1)^d)$ such that $0\le\eta\le1$ and such that
$\eta(s)=1$ if and only if $s\in F$.
\end{lemma}

We say that a sequence $(a_n)_{n\in\Z^d}$ of complex numbers is \emph{rapidly decreasing} if $\sup_{n\in\Z^d}(1+|n|)^k|a_n|<\infty$ for every $k\in\N$.

\begin{proposition}
\label{prop:exact-gap}
Let $F\subset\R^d$ be compact and contained in an open cube $(0,1)^d-\tau$ for some $\tau\in\R^d$.  Then there
are rapidly decreasing
sequences $(a_n)_{n\in\Z^d}$ and $(b_n)_{n\in\Z^d}$ for which
\[
 \Phi_F(\omega)
 =\sum_{n\in\Z^d}a_ne^{-2\pi in\cdot\omega}
 +\sum_{n\in\Z^d}b_ne^{-2\pi i(n+\vartheta)\cdot\omega}
\]
converges uniformly on $\R^d$ and satisfies
\[
 \{\omega:\Phi_F(\omega)=0\}=F.
\]
Moreover, the coefficients may be chosen so that $b_{\mathbf0}\neq0$.
For $F=\varnothing$, one may take $\Phi_F(\omega)=1+\frac{1}{2}e^{-2\pi i \vartheta\cdot\omega}$.
\end{proposition}

\begin{proof} Translate $F$ by a vector $\tau$ so that
\[
 F_0:=F+\tau\subset(0,1)^d.
\]
Let $\eta$ be given by Lemma~\ref{lem:smooth-peak}, and define $g_0:[0,1)^d\to\C$ via
\[
 g_0(s)=-\eta(s)e^{-2\pi i\vartheta\cdot s}.
\]
Since $\eta$ vanishes in a neighborhood of the boundary of the cube $[0,1]^d$, the $\Z^d$-periodic extension of
$g_0$ is smooth.  Hence, $g_0$ can be written as
\[
 g_0(s)=\sum_{n\in\Z^d}c_ne^{-2\pi i n\cdot s}
\]
with rapidly decreasing coefficients $(c_n)$.  For $\omega=m+s$ with $m\in\Z^d$ and $s\in[0,1)^d$,
set
\[
 \Phi_0(\omega)=g_0(\omega)+e^{-2\pi i\vartheta\cdot\omega}=e^{-2\pi i\vartheta\cdot s}
   \bigl(e^{-2\pi i\vartheta\cdot m}-\eta(s)\bigr).
\]
Because $0\le\eta\le1$, the right-hand side can vanish only if
$\eta(s)=1$ and $e^{-2\pi i\vartheta\cdot m}=1$.  The rational
independence of $1,\vartheta_1,\ldots,\vartheta_d$ forces $m=0$, and thus
$\{\omega\in\R^d:\Phi_0(\omega)=0\}=F+\tau$.  Finally, put
\[
 \Phi_F(\omega)=\Phi_0(\omega+\tau).
\]
Then $\Phi_F$ is of the asserted form with the coefficients $a_n=c_ne^{-2\pi in\cdot\tau}$, $b_{\mathbf0}=e^{-2\pi i\vartheta\cdot\tau}$ and $b_n=0$ for $n\neq\mathbf0$. In particular, we have $b_{\mathbf0}\neq0$.
\end{proof}

Let $\mathcal E$ denote the class of all series of the form
\begin{equation}
\label{eq:psiomega}
 \Psi(\omega)
 =\sum_{n\in\Z^d}a_ne^{-2\pi in\cdot\omega}
 +\sum_{n\in\Z^d}b_ne^{-2\pi i(n+\vartheta)\cdot\omega}
\end{equation}
with rapidly decreasing coefficients $(a_n)$ and $(b_n)$.  Since
$\Z^d\cap(\Z^d+\vartheta)=\emptyset$ we may write this sum
as a single absolutely convergent sum
\[
 \Psi(\omega)=\sum_{\lambda\in\Lambda}
 c_\lambda^\Psi e^{-2\pi i\lambda\cdot\omega},
\]
where $c_\lambda^\Psi=a_n$ for $\lambda=n$ and
$c_\lambda^\Psi=b_n$ for $\lambda=n+\vartheta$.

The class $\mathcal E$ is stable under
smooth periodic multipliers, and such a multiplier can be arranged to
agree with a chosen exponential on one small cube.

\begin{lemma}
\label{lem:multipliers}
The following statements hold.
\begin{enumerate}
\item If $\Psi\in\mathcal E$ and $r\in\mathcal{C}^\infty(\mathbb{T}^d)$,
then $r\Psi\in\mathcal E$.

\item Let $J=\tau+[-\ell/2,\ell/2]^d$ with $\tau\in\R^d$ and $0<\ell<1$. For every
      $y\in\R^d$, there exists a smooth
      $\Z^d$-periodic function $r_y$ such that
      $r_y(\omega)=e^{2\pi i y\cdot\omega}$
      on a neighborhood of $J$.
\end{enumerate}
\end{lemma}

\begin{proof}
(1) A smooth $\Z^d$-periodic function $r$ has rapidly decreasing Fourier
coefficients $(\widehat{r}(n))$ supported on $\Z^d$, and we may write 
$$
r(\omega)=\sum_{n\in\Z^d}\widehat{r}(n)e^{2i\pi n\cdot\omega}.
$$
Therefore if $\Psi$ is of the form \eqref{eq:psiomega} then, with $c_n=\widehat{r}(-n)$,
$$
r(\omega) \Psi(\omega)=\sum_{n\in\Z^d}(a*c)(n)e^{-2\pi in\cdot\omega}
 +\sum_{n\in\Z^d}(b*c)(n)e^{-2\pi i(n+\vartheta)\cdot\omega}
$$
and, if $a,b$ is rapidly decreasing, so are $a*c$ and $b*c$. Multiplication by $\mathcal{C}^\infty(\mathbb{T}^d)$ therefore preserves $\mathcal E$.

\smallskip

(2) Choose $\chi\in C_c^\infty(\tau+(-1/2,1/2)^d)$ with
$\chi=1$ on a neighborhood of $J$, and set
\[
 r_y(\omega)=\sum_{n\in\Z^d}
 \chi(\omega-n)e^{2\pi i y\cdot(\omega-n)}.
\]
The function $r_y$ is smooth and $\Z^d$-periodic.  Near $J$, only the summand $n=0$
is nonzero, giving the asserted identity.
\end{proof}

\section{Local realization near the delay axis}
\label{sec:local}

In this section, we are going to prove Theorem \ref{thm:main}. We start by reducing it
to a more localized statement.

\subsection{Symplectic reduction}

We first show that there exists a symplectic transform that maps a small ball around any
nonzero point $a$ into a small cube centered at a point $s_0e_1$ of the delay
axis. By the symplectic covariance \eqref{eq:full-symplectic-covariance}, this
reduces Theorem~\ref{thm:main} to holes located near the delay axis.

\begin{lemma}\label{lem:sumpreduc}
Let $Q=s_0 e_1+[-\rho,\rho]^{2d}\subset\mathbb R^{2d}$ with $s_0>0$ and
$0<\rho<s_0$. Further, let \(a\in\mathbb R^{2d}\setminus\{0\}\). There exists a linear symplectic transform 
\(S\in \operatorname{Sp}(2d,\mathbb R)\) and $r_a>0$ such that,
\[
S\bigl(B(a,r_a)\bigr)\subset Q.
\]
\end{lemma}

\begin{proof}
We identify \(\mathbb R^{2d}\) with \(\mathbb C^d\), endowed with its standard symplectic structure. Since
$O(2d)\cap \operatorname{Sp}(2d,\mathbb R)$
acts transitively on the unit sphere of \(\mathbb R^{2d}\), there exists
\[
U\in O(2d)\cap \operatorname{Sp}(2d,\mathbb R)
\]
such that $Ua=|a|e_1$.

Set $\lambda=s_0/|a|$ and consider the matrix
$D_\lambda=\operatorname{diag}(\lambda,1,\ldots,1,\lambda^{-1},1,\ldots,1)$
with respect to the symplectic coordinates $(q_1,\ldots,q_d,p_1,\ldots,p_d)$.
The matrix $D_\lambda$ is symplectic. Hence
$S=D_\lambda U\in\operatorname{Sp}(2d,\mathbb R)$, and
$Sa=D_\lambda(|a|e_1)=s_0e_1$. Since $U$ is orthogonal, we have
$\|S\|=\|D_\lambda\|=\max\{\lambda,\lambda^{-1}\}=\max\{s_0/|a|,|a|/s_0\}$.
Therefore, if $x\in B(a,r)$, then
$|Sx-s_0e_1|=|S(x-a)|\leq\|S\|\,|x-a|<\|S\|r$, that is,
$S(B(a,r))\subset B(s_0e_1,\|S\|r)$. Since
$B(s_0e_1,\rho)\subset s_0e_1+[-\rho,\rho]^{2d}=Q$, it is enough to require
$\|S\|r\leq\rho$, that is, $r\leq\rho\min\{|a|/s_0,s_0/|a|\}$. Thus,
setting $r_a=\rho\min\{|a|/s_0,s_0/|a|\}$, we obtain
$S(B(a,r_a))\subset Q$.
\end{proof}

The proof of the local realization is given at the end of this section.
We first show how it implies Theorem~\ref{thm:main}.

\begin{proof}[Proof of Theorem~\ref{thm:main} assuming the local realization below]
Fix $a\in\R^{2d}\setminus\{0\}$, and let $\Omega\subset\R^{2d}$ be
bounded, nonempty, and regular open. Put
\[
R_\Omega=\sup_{u\in\Omega}|u|
\]
and choose $b_0>0$ with $b_0R_\Omega<r_a$, where $r_a>0$ and $S\in\Sp(2d,\R)$ are given by
Lemma~\ref{lem:sumpreduc}. For every $0<b<b_0$, we have $a+b\Omega\subset B(a,r_a)$, and hence the
map $S$ satisfies
\[
S(a+b\Omega)\subset Q=s_0e_1+[-\rho,\rho]^{2d}.
\]
Using that $S(a+b\Omega)$ is regular open, an application of
the local realization below gives $f\in L^2(\R^d)$ with
\[
\Zopen(f)=S(a+b\Omega)\cup(-S(a+b\Omega)).
\]
The symplectic covariance
$\Zopen(\mathcal U_{S^{-1}}f)=(a+b\Omega)\cup(-(a+b\Omega))$
yields the statement.
\end{proof}

\begin{corollary}
\label{cor:quantitative}
For every $d\in\N$, there exists a constant $\kappa_d>0$, depending
only on $d$, with the following property. Let $a\in\R^{2d}\setminus\{0\}$,
let $\Omega\subset\R^{2d}$ be bounded, nonempty, and regular open, and set
$R_\Omega:=\sup_{u\in\Omega}|u|$. Then, for every
$0<b<\kappa_d\min\{|a|,|a|^{-1}\}/R_\Omega$, there exists
$f\in L^2(\R^d)$ such that
\[
\Zopen(f)=(a+b\Omega)\cup(-a-b\Omega).
\]
\end{corollary}

\begin{proof}
Take $\rho$ and $s_0$ from \eqref{eq:parameters}, and set
$\kappa_d=\rho\min\{s_0,s_0^{-1}\}$.
Then $bR_\Omega<r_a$ and therefore $a+b\Omega\subset B(a,r_a)$.
The preceding proof yields the statement.
\end{proof}

\subsection{The functions \texorpdfstring{$W$}{W} and \texorpdfstring{$V^\Psi$}{V}}

For $\lambda\in\Lambda=\Z^d\cup(\vartheta+\Z^d)$, write
$\lambda=n$ or $\lambda=n+\vartheta$ and put
\[
w_\lambda=(1+|n|)^{-(d+1)}.
\]

Further, let $\mathcal Q=[0,1]^d$ be the unit cube. If
$L<\operatorname{sep}_\infty(\Lambda)$, then the cubes
$\{L\mathcal Q+\lambda:\lambda\in\Lambda\}$ are pairwise disjoint.
In particular, the function
\[
W(t)
:=\sum_{\lambda\in\Lambda}w_\lambda\1_{L\qq}(t-\lambda)
=\sum_{\lambda\in\Lambda}w_\lambda\tau_\lambda\1_{L\qq}(t)
\]
belongs to $L^2(\R^d)$. To the function
\[
\Psi(\omega)
=\sum_{\lambda\in\Lambda}
c_\lambda^\Psi e^{-2\pi i\lambda\cdot\omega}
\in\mathcal E,
\]
we associate
\[
v_\lambda^\Psi=\frac{c_\lambda^\Psi}{w_\lambda},
\qquad
V^\Psi(t)
=\sum_{\lambda\in\Lambda}v_\lambda^\Psi\1_{L\qq}(t-\lambda)
=\sum_{\lambda\in\Lambda}v_\lambda^\Psi
\tau_\lambda\1_{L\qq}(t).
\]
The rapid decrease of the coefficients $(c_\lambda^\Psi)_{\lambda\in\Lambda}$ and the polynomial growth of
$(w_\lambda^{-1})_{\lambda\in\Lambda}$ show that
$(v_\lambda^\Psi)_{\lambda\in\Lambda}\in\ell^2(\Lambda)$.
Thus $V^\Psi\in L^2(\R^d)$.

\begin{lemma}
\label{lem:train}
Let $L>4\rho>0$, put $\sigma=\sep(\Lambda)$, and assume
\begin{equation}
3L+7\rho<s_0,
\qquad
2s_0+2L+4\rho<\sigma.
\label{eq:conditions0}
\end{equation}
Then, for
\[
\alpha,\beta\in
T:=s_0e_1+[-(L+3\rho),L+3\rho]^d
\]
and $\Psi,\Phi\in\mathcal E$, the following hold:
\begin{enumerate}[label=\textup{(\roman*)}]
\item
$A(W)(x,\omega)=0$ for
$x\in\pm s_0e_1+[-\rho,\rho]^d$.

\item
$A(\tau_\alpha V^\Psi,\tau_\beta V^\Phi)(x,\omega)=0$
for $x\in\pm s_0e_1+[-\rho,\rho]^d$.

\item
$A(W,\tau_\alpha V^\Psi)(x,\omega)=0$
for $x\in s_0e_1+[-\rho,\rho]^d$.

\item
For $x\in-s_0e_1+[-\rho,\rho]^d$,
\[
A(W,\tau_\alpha V^\Psi)(x,\omega)
=e^{-\pi i\omega\cdot\alpha}
\overline{\Psi(-\omega)}
A(\1_{L\qq})(x+\alpha,\omega).
\]

\item
$A(\tau_\alpha V^\Psi,W)(x,\omega)=0$
for $x\in-s_0e_1+[-\rho,\rho]^d$.

\item
For $x\in s_0e_1+[-\rho,\rho]^d$,
\[
A(\tau_\alpha V^\Psi,W)(x,\omega)
=e^{-\pi i\omega\cdot\alpha}
\Psi(\omega)A(\1_{L\qq})(x-\alpha,\omega).
\]
\end{enumerate}
\end{lemma}

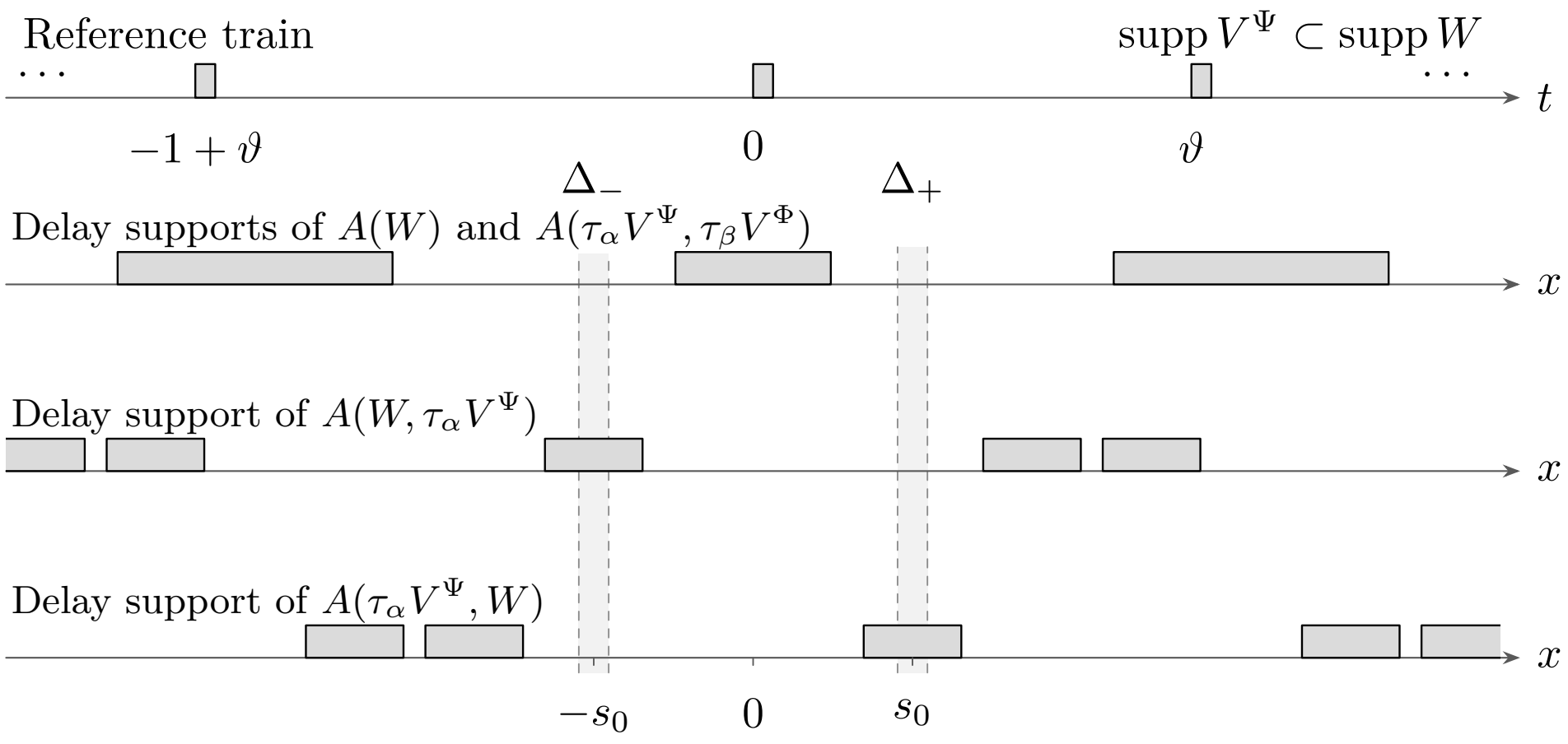
\begin{figure}[htbp]
\centering

\tikzset{
  paperfigure/.style={
    font=\small,
    line cap=round,
    line join=round
  },
  figaxis/.style={
    draw=black!65,
    line width=.4pt,
    -{Stealth[length=1.4mm,width=1.05mm]}
  },
  figblock/.style={
    draw=black,
    fill=black!14,
    line width=.45pt
  },
  figguide/.style={
    draw=black!45,
    line width=.35pt,
    dash pattern=on 2pt off 2pt
  }
}

\begin{tikzpicture}[paperfigure,x=8cm,y=1cm]
  \node[anchor=west] at (-.75,5.02) {Reference train};
  \node[anchor=east] at (.75,5.02)
    {$\supp V^\Psi\subset\supp W$};

  \draw[figaxis] (-.75,4.5)--(.77,4.5)
    node[right] {$t$};

  \foreach \c in {-.56,0,.44}{
    \path[figblock]
      (\c,4.5) rectangle ({\c+.02},4.77);
  }
  \foreach \c/\lab in {
    -.56/{-1+\vartheta},
    0/0,
    .44/\vartheta
  }{
    \node[below=4pt] at (\c,4.5) {$\lab$};
  }
  \node at (-.71,4.68) {$\cdots$};
  \node at (.70,4.68) {$\cdots$};

  \foreach \c in {-.16,.16}{
    \fill[black!5]
      ({\c-.015},-.12) rectangle ({\c+.015},3.30);
    \draw[figguide]
      ({\c-.015},-.12)--({\c-.015},3.30);
    \draw[figguide]
      ({\c+.015},-.12)--({\c+.015},3.30);
  }
  \node at (-.16,3.82) {$\Delta_-$};
  \node at (.16,3.82) {$\Delta_+$};

  \node[
    anchor=west,
    fill=white,
    inner sep=1pt,
    font=\footnotesize
  ] at (-.75,3.46)
    {Delay supports of $A(W)$ and
     $A(\tau_\alpha V^\Psi,\tau_\beta V^\Phi)$};

  \draw[figaxis] (-.75,3)--(.77,3)
    node[right] {$x$};
  \path[figblock] (-.638,3) rectangle (-.362,3.26);
  \path[figblock] (-.078,3) rectangle (.078,3.26);
  \path[figblock] (.362,3) rectangle (.638,3.26);

  \node[
    anchor=west,
    fill=white,
    inner sep=1pt,
    font=\footnotesize
  ] at (-.75,1.96)
    {Delay support of $A(W,\tau_\alpha V^\Psi)$};

  \draw[figaxis] (-.75,1.5)--(.77,1.5)
    node[right] {$x$};
  \begin{scope}
    \clip (-.75,1.49) rectangle (.75,1.77);
    \foreach \c in {-.56,-.44,0,.44,.56}{
      \path[figblock]
        ({\c-.16-.049},1.5)
        rectangle ({\c-.16+.049},1.76);
    }
  \end{scope}

  \node[
    anchor=west,
    fill=white,
    inner sep=1pt,
    font=\footnotesize
  ] at (-.75,.46)
    {Delay support of $A(\tau_\alpha V^\Psi,W)$};

  \draw[figaxis] (-.75,0)--(.77,0)
    node[right] {$x$};
  \begin{scope}
    \clip (-.75,-.01) rectangle (.75,.27);
    \foreach \c in {-.56,-.44,0,.44,.56}{
      \path[figblock]
        ({\c+.16-.049},0)
        rectangle ({\c+.16+.049},.26);
    }
  \end{scope}

  \foreach \c/\lab in {-.16/{-s_0},0/0,.16/s_0}{
    \draw[draw=black!65,line width=.35pt]
      (\c,0)--(\c,-.06);
    \node[below=4pt] at (\c,-.06) {$\lab$};
  }
\end{tikzpicture}

\caption{Support separation in dimension one. The vertical strips mark
the sets $\Delta_\pm=\pm s_0+[-\rho,\rho]$; the horizontal blocks bound the relevant supports.
The pure terms miss both strips. The mixed terms can meet only
$\Delta_-$ and $\Delta_+$, respectively.
Spacings and widths are schematic.}
\label{fig:train-separation}
\end{figure}

\begin{proof}
Let $r,s\in\R^d$. By Lemma~\ref{lem:basic-facts}\ref{fact:translation}, we have
\begin{equation}\label{eq:cube-pair}
A(\tau_r\1_{L\qq},\tau_s\1_{L\qq})(x,\omega)=e^{-\pi i\omega\cdot(r+s)}A(\1_{L\qq})(x-(r-s),\omega),
\end{equation}
and by Lemma~\ref{lem:basic-facts}\ref{fact:expra1q}, the right-hand side of \eqref{eq:cube-pair} vanishes for $x\notin(r-s)+[-L,L]^d$. Moreover, since the series defining $W$ and $V^\Psi$ converge in $L^2(\R^d)$ and have absolutely summable coefficients, Lemma~\ref{lem:basic-facts}\ref{fact:continuity} shows that all series of ambiguity functions below converge absolutely and uniformly on $\R^{2d}$.

We first prove (i) and (ii). To do so, let $x\in\pm s_0e_1+[-\rho,\rho]^d$, so that $s_0-\rho\le\|x\|_\infty\le s_0+\rho$. By \eqref{eq:cube-pair}, we have
$$
A(\tau_\alpha V^\Psi,\tau_\beta V^\Phi)(x,\omega)=\sum_{r,s\in\Lambda}v_r^\Psi\overline{v_s^\Phi}\,e^{-\pi i\omega\cdot(\alpha+\beta+r+s)}A(\1_{L\qq})(x-(\alpha-\beta)-(r-s),\omega).
$$
Suppose that the term with indices $r,s\in\Lambda$ does not vanish, and set $\nu=r-s$. Since $\alpha-\beta\in[-(2L+6\rho),2L+6\rho]^d$, it follows that $\|x-\nu\|_\infty\le3L+6\rho$. If $\nu=0$, then $s_0-\rho\le3L+6\rho$, which contradicts \eqref{eq:conditions0}. If $\nu\neq0$, then $\|\nu\|_\infty\ge\sigma$, and therefore
$$
\sigma\le\|x\|_\infty+3L+6\rho\le s_0+3L+7\rho<2s_0<\sigma,
$$
which gives a contradiction. Hence, every term vanishes, which proves (ii). The same argument with $\alpha=\beta=0$ and with $v_r^\Psi,v_s^\Phi$ replaced by $w_r,w_s$ proves (i).

Next, we prove (iii) and (iv). By \eqref{eq:cube-pair}, we have
\begin{equation}\label{eq:mixed-sum}
A(W,\tau_\alpha V^\Psi)(x,\omega)=e^{-\pi i\omega\cdot\alpha}\sum_{r,s\in\Lambda}w_r\overline{v_s^\Psi}\,e^{-\pi i\omega\cdot(r+s)}A(\1_{L\qq})(x+\alpha-(r-s),\omega),
\end{equation}
and the term with indices $r,s\in\Lambda$ vanishes unless $\nu=r-s$ satisfies $\nu\in x+\alpha+[-L,L]^d$.

If $x\in s_0e_1+[-\rho,\rho]^d$, then $x+\alpha\in2s_0e_1+[-(L+4\rho),L+4\rho]^d$, and therefore every $\nu\in x+\alpha+[-L,L]^d$ satisfies
$$
\nu_1\ge2s_0-(2L+4\rho)>0,\qquad\|\nu\|_\infty\le2s_0+2L+4\rho<\sigma.
$$
Since no $\nu\in\Lambda-\Lambda$ satisfies both inequalities, every term in \eqref{eq:mixed-sum} vanishes, which proves (iii).

If $x\in-s_0e_1+[-\rho,\rho]^d$, then $x+\alpha\in[-(L+4\rho),L+4\rho]^d$, and therefore every $\nu\in(\Lambda-\Lambda)\cap(x+\alpha+[-L,L]^d)$ satisfies $\|\nu\|_\infty\le2L+4\rho<\sigma$, that is, $\nu=0$. Hence, only the terms with $r=s$ remain in \eqref{eq:mixed-sum}, and since $w_r\overline{v_r^\Psi}=\overline{c_r^\Psi}$, we obtain
$$
\begin{aligned}
A(W,\tau_\alpha V^\Psi)(x,\omega)&=e^{-\pi i\omega\cdot\alpha}\sum_{r\in\Lambda}\overline{c_r^\Psi}e^{-2\pi ir\cdot\omega}A(\1_{L\qq})(x+\alpha,\omega)\\
&=e^{-\pi i\omega\cdot\alpha}\overline{\Psi(-\omega)}A(\1_{L\qq})(x+\alpha,\omega),
\end{aligned}
$$
which proves (iv).

Finally, Lemma~\ref{lem:basic-facts}\ref{fact:sym} gives
$$
A(\tau_\alpha V^\Psi,W)(x,\omega)=\overline{A(W,\tau_\alpha V^\Psi)(-x,-\omega)},\qquad\overline{A(\1_{L\qq})(-x+\alpha,-\omega)}=A(\1_{L\qq})(x-\alpha,\omega).
$$
Therefore, (v) follows from (iii), and (vi) follows from (iv).
\end{proof}

For the remainder of this section, put $\sigma=\sep(\Lambda)$ and fix the constants
\begin{equation}
s_0=\frac{\sigma}{4},
\qquad
L=\frac{\sigma}{100},
\qquad
\rho=\frac{\sigma}{1000}.
\label{eq:parameters}
\end{equation}
These choices satisfy the hypotheses of Lemma~\ref{lem:train}. Moreover, we write
\[
\Delta_+=s_0e_1+[-\rho,\rho]^d,
\qquad
\Delta_-=-s_0e_1+[-\rho,\rho]^d,
\]
and $\Delta=\Delta_+\cup\Delta_-$.

\subsection{Separating one exterior point}
\label{subsec:separating}

We will need the following simple computation.

\begin{lemma}
\label{lem:one-coordinate-difference}
Let $0<h<L/2$ and define the function $\mathcal B_{L,h}:\R^2\to\C$ via
\[
 \mathcal B_{L,h}(x,\omega)
 =A(\1_{[0,L]})(x,\omega)
 -2\cos(\pi \omega h)A(\1_{[0,L]})(x-h,\omega)
 +A(\1_{[0,L]})(x-2h,\omega).
\]
Then, for every $\omega\in\R$, we have
\[
 \supp_x\mathcal B_{L,h}(\cdot,\omega)
 \subset[-L,-L+2h]\cup[0,2h]\cup[L,L+2h],
\]
and
\[
 \mathcal B_{L,h}(L+h,\omega)=he^{-\pi iL\omega }\sinc(\pi\omega h).
\]
\end{lemma}

\begin{proof}
First $\supp_xA(\1_{[0,L]})(\cdot,\omega)=[-L,L]$ so that 
$\supp_x\mathcal B_{L,h}(\cdot,\omega)\subset[-L,L+2h]$.

Now let $\omega\neq0$. By Lemma~\ref{lem:basic-facts}\ref{fact:expra1q}, we have
$$
A(\1_{[0,L]})(x,\omega)=\frac{e^{-\pi iL\omega}}{\pi\omega}\begin{cases}\sin(\pi(L+x)\omega),&x\in[-L,0],\\ \sin(\pi(L-x)\omega),&x\in[0,L].\end{cases}
$$
Since every function of the form $u(x)=\sin(a+bx)$ with $a\in\R$ and $b\in\{\pi\omega,-\pi\omega\}$ satisfies $u(x)-2\cos(\pi\omega h)u(x-h)+u(x-2h)=0$ for all $x\in\R$, it follows that $\mathcal B_{L,h}(x,\omega)=0$ if $x,x-h,x-2h\in[-L,0]$ or $x,x-h,x-2h\in[0,L]$, that is, if $x\in[-L+2h,0]\cup[2h,L]$. This proves the inclusion for $\omega\neq0$, and the case $\omega=0$ follows from the continuity of $\mathcal B_{L,h}$.

Finally, since $A(\1_{[0,L]})(L+h,\omega)=A(\1_{[0,L]})(L,\omega)=0$, Lemma~\ref{lem:basic-facts}\ref{fact:expra1q} gives $\mathcal B_{L,h}(L+h,\omega)=A(\1_{[0,L]})(L-h,\omega)=he^{-\pi iL\omega}\sinc(\pi\omega h)$.
\end{proof}

For $\mathbf h=(h_1,\ldots,h_d)$ define
\[
 \mathcal B_{L,\mathbf h}^{(d)}(x,\omega)
 :=\prod_{j=1}^d\mathcal B_{L,h_j}(x_j,\omega_j).
\]
Its delay support lies in the union of $3^d$ boxes obtained by independently
choosing one of the three intervals $(-L,-L+2h_j)$, $(0,2h_j)$ or $(L,L+2h_j)$ in each coordinate. 
A suitable translation will ensure that only one of those boxes intersects $s_0e_1+[-\rho,\rho]^d$.

\begin{lemma}
\label{lem:local-packet}
Let $\Omega\subset \Delta_+ \times [-\rho,\rho]^d$ be a regular open set and let $K=\overline{\Omega}$.
Then for every
\[
 z=(y,\eta)\in(\Delta_+\times\R^d)\setminus K,
\]
there is a sum $R_z$ of at most $3^d$ functions belonging to
$$\{ \tau_\alpha V^\Psi : \alpha\in s_0e_1+[-(L+3\rho),L+\rho]^d, \, \Psi\in\mathcal E \}$$ such that
\[
 A(R_z,W)(x,\omega)=0\quad\text{for }(x,\omega)\in K,
 \qquad A(R_z,W)(y,\eta)\neq0.
\]
\end{lemma}

\begin{proof}
The construction depends on whether the target frequency $\eta$ belongs to $[-\rho,\rho]^d$ or not.

\medskip

\noindent{\bf Case 1.} {\em $\eta\notin [-\rho,\rho]^d$.}

\medskip

We first chose $\beta=(\beta_1,\ldots,\beta_d)$ with $\beta_j\in[-L,L]$ such that
\[
A(\1_{[0,L]})(\beta_j,\eta_j)=e^{-\pi iL\eta_j}(L-|\beta_j|)\sinc\bigl(\pi(L-|\beta_j|)\eta_j\bigr) \not=0.
\]
We then set $\alpha=y-\beta$
so that $A(\1_{L\qq})(y-\alpha,\eta)\neq0$. Note that $\alpha\in s_0e_1+[-L-\rho,L+\rho]^d$.

Now put
\[
 F=\{\omega\in\R^d:\exists x\in\alpha+[-L,L]^d
       \text{ with }(x,\omega)\in K\}.
\]
The set $F$ is compact and lies in $[-\rho,\rho]^d$, so $\eta\notin F$. 
We can thus apply Proposition~\ref{prop:exact-gap} which provides us with
$\Psi=\Phi_F$ such that $\Psi(\omega)=0$ if and only if $\omega\in F$. In particular, we have $\Psi(\eta)\not=0$.
Set $R_z=\tau_\alpha V^\Psi$.  

Now, for $x\in \Delta_+$, Lemma~\ref{lem:train}, implies that
$$
A(\tau_\alpha V^\Psi,W)(x,\omega)=e^{-\pi i \omega \cdot \alpha}\Phi_F(\omega)A(\1_{L\qq})(x-\alpha,\omega).
$$
If $x\notin \alpha+[-L,L]^d$, then $A(\tau_\alpha V^\Psi,W)(x,\omega)=0$
and if $x\in \alpha+[-L,L]^d$ and $(x,\omega)\in K$ then $\omega\in F$ thus $\Psi(\omega)=0$
and $A(\tau_\alpha V^\Psi,W)(x,\omega)=0$. This shows that $A(\tau_\alpha V^\Psi,W)=0$ on $K$.

On the other hand, $y\in \alpha+[-L,L]^d$ and $\eta \not\in [-\rho,\rho]^d$, thus $\eta\notin F$ and $\Psi(\eta)\not=0$
and $A(\1_{L\qq})(y-\alpha,\eta)\not=0$ by construction of $\alpha$. Hence, it follows that
$$
A(\tau_\alpha V^\Psi,W)(y,\eta)=e^{-\pi i\eta \cdot \alpha}\Phi_F(\eta)A(\1_{L\qq})(y-\alpha,\eta)\not=0.
$$

\medskip

\noindent{\bf Case 2.} {\sl $\eta\in [-\rho,\rho]^d$.}

\medskip

Consider the set
\[
 S_\eta=\{x\in\R^d:(x,\eta)\in K\},
\]
which is compact and does not contain $y$.  Choose numbers
$0<h_j<\rho$ so small that $h_j<1/|\eta_j|$ when $\eta_j\not=0$ and such that
\begin{equation}
 y+\prod_{j=1}^d[-h_j,h_j]
 \quad\text{is disjoint from }S_\eta.
 \label{eq:selected-box-disjoint}
\end{equation}
Note that, whenever $\eta_j\neq0$, $\sin(\pi\eta_jh_j)\neq0$.  Put
\[
 \mathbf h=(h_1,\ldots,h_d),
 \qquad \xi=y-L\one-\mathbf h,
\]
where $\one=(1,\ldots,1)\in\R^d$.
Writing $\odot$ for componentwise multiplication on $\R^d$, we observe that
every vector $\xi+k\odot\mathbf h$, $k\in\{0,1,2\}^d$, lies in 
$$
s_0e_1+[-\rho,\rho]^d-L\one-[-\rho,\rho]^d\subset s_0e_1+[-(L+3\rho),L+\rho]^d.
$$

Define
\[
 F=\left\{\omega:\exists u\in  y+\prod_{j=1}^d[-h_j,h_j] \text{ with }(u,\omega)\in K\right\}.
\]
Note that $F$ is compact, $F\subset [-\rho,\rho]^d$ and $\eta\notin F$ by \eqref{eq:selected-box-disjoint}.
Let $\Phi_F$ be the function associated to $F$ in Proposition \ref{prop:exact-gap}.

For each $j$, use Lemma~\ref{lem:multipliers} to choose a smooth periodic
$r_j$ satisfying
\[
 r_j(\omega)=e^{2\pi i h_j\omega_j},
 \qquad\omega\in [-\rho,\rho]^d.
\]
For $r\in\C$, set
\[
 d_0(r)=1,
 \qquad d_1(r)=-(1+r),
 \qquad d_2(r)=r,
\]
and for $k\in\{0,1,2\}^d$ set
\[
 \Psi_k(\omega) =\Phi_F(\omega)\prod_{j=1}^d d_{k_j}\bigl(r_j(\omega)\bigr),
 \qquad
 R_z=\sum_{k\in\{0,1,2\}^d}\tau_{\xi+k\odot\mathbf h}V^{\Psi_k}.
\]
By Lemma~\ref{lem:multipliers}, we have $\Psi_k\in\mathcal E$ for every $k\in\{0,1,2\}^d$.

We proceed by showing that
\begin{equation}
 A(R_z,W)(x,\omega)
 =e^{-\pi i\omega\cdot\xi}\Phi_F(\omega)
   \mathcal B_{L,\mathbf h}^{(d)}(x-\xi,\omega),
 \quad (x,\omega)\in\Delta_+\times[-\rho,\rho]^d.
 \label{eq:localized-product}
\end{equation}
To do so, let $(x,\omega)\in\Delta_+\times[-\rho,\rho]^d$ and $k\in\{0,1,2\}^d$. Lemma~\ref{lem:train}(vi) and Lemma~\ref{lem:basic-facts}\ref{fact:expra1q} give
$$
A(\tau_{\xi+k\odot\mathbf h}V^{\Psi_k},W)(x,\omega)=e^{-\pi i\omega\cdot\xi}\Phi_F(\omega)\prod_{j=1}^de^{-\pi i\omega_jk_jh_j}d_{k_j}\bigl(e^{2\pi ih_j\omega_j}\bigr)A(\1_{[0,L]})(x_j-\xi_j-k_jh_j,\omega_j).
$$
Summing over $k\in\{0,1,2\}^d$ and using the identities $e^{-\pi i\omega_jh_j}d_1(e^{2\pi ih_j\omega_j})=-2\cos(\pi\omega_jh_j)$ and $e^{-2\pi i\omega_jh_j}d_2(e^{2\pi ih_j\omega_j})=1$, we obtain \eqref{eq:localized-product}.

Next, we show that $A(R_z,W)=0$ on $K$. To do so, let $(x,\omega)\in K$ and suppose that $\mathcal B_{L,\mathbf h}^{(d)}(x-\xi,\omega)\neq0$. By Lemma~\ref{lem:one-coordinate-difference} and the definition of $\xi$, for every $j\in\{1,\dots,d\}$ we have
$$
x_j-y_j\in[-2L-h_j,-2L+h_j]\cup[-L-h_j,-L+h_j]\cup[-h_j,h_j].
$$
Since $x,y\in\Delta_+$ implies $|x_j-y_j|\le2\rho$, and since $-L+h_j<-3\rho$, it follows that $x_j-y_j\in[-h_j,h_j]$. Hence, $x\in y+\prod_{j=1}^d[-h_j,h_j]$, and therefore $\omega\in F$ and $\Phi_F(\omega)=0$. Consequently, \eqref{eq:localized-product} shows that $A(R_z,W)(x,\omega)=0$ for all $(x,\omega)\in K$.

Finally, since $y-\xi=L\one+\mathbf h$, it follows from \eqref{eq:localized-product} and Lemma~\ref{lem:one-coordinate-difference} that
$$
A(R_z,W)(y,\eta)=e^{-\pi i\eta\cdot\xi}\Phi_F(\eta)\prod_{j=1}^de^{-\pi i\eta_jL}h_j\sinc(\pi\eta_jh_j).
$$
The right-hand side is nonzero, since $\eta\notin F$ and $|\eta_jh_j|<1$ for every $j\in\{1,\dots,d\}$.
\end{proof}

\subsection{Support near the selected delays}
\label{subsec:inner-support}

We now apply Lemma~\ref{lem:local-packet} to a countable dense subset of
the complement of $K$ in $\Delta_+\times\R^d$.
Choose a countable dense set
\[
 Z_+=\{z_n:n\ge1\}  \subset\bigl(\Delta_+\times\R^d)\setminus K.
\]
For each $n$, choose $R_n=R_{z_n}$ from Lemma~\ref{lem:local-packet}.  By
Lemma~\ref{lem:train}, for arbitrary real $t_1,\ldots,t_N$, we have
\[
 A\left(W+\sum_{j=1}^Nt_jR_j\right) =\sum_{j=1}^Nt_jA(R_j,W)
  \qquad\text{on }\Delta_+\times\R^d.
\]
Thus every ambiguity function of the form $A(W+\sum_{j=1}^Nt_jR_j)$ vanishes on $K$ and $-K$.  At $z_n$, adding $tR_n$ gives a polynomial with nonzero
linear coefficient $A(R_n,W)(z_n)$.  Lemma~\ref{lem:diagonal} yields real
values $t_n$ such that the series
\begin{equation}
 H=W+\sum_{n=1}^\infty t_nR_n
 \label{eq:H}
\end{equation}
converges in $L^2(\R^d)$, and satisfies the properties
$$
 A(H)=0\quad\text{on }K \cup -K,
 \qquad A(H)(z_n)\neq0,\quad n\ge1.
$$
Every $z_n$ belongs to the support of $A(H)$, which is a closed set.  Since
$Z_+$ is dense in
$(\Delta_+\times\R^d)\setminus K$, it follows that the support of $A(H)$
contains the closure of $(\Delta_+\times\R^d)\setminus K$.  In particular, we have
\[
 \overline{(\Delta_+\times\R^d)\setminus K}
 \supset (\Delta_+\times\R^d)\setminus K,
\]
which gives
\begin{equation}
 (\Delta_+\times\R^d)\setminus K\subset\supp A(H),
 \qquad
 (\Delta_-\times\R^d)\setminus -K\subset\supp A(H).
 \label{eq:inner-support}
\end{equation}

\begin{figure}[htbp]
\centering
\begin{tikzpicture}[paperfigure,x=1.65cm,y=1.35cm]
  \foreach \l/\r in {-2.7/-1.1,1.1/2.7}{
    \fill[black!5] (\l,-1.3) rectangle (\r,1.3);
    \draw[figguide] (\l,-1.3)--(\l,1.3);
    \draw[figguide] (\r,-1.3)--(\r,1.3);
  }
  \draw[figaxis] (-3.6,0)--(3.6,0) node[right] {$x$};
  \draw[figaxis] (0,-1.5)--(0,1.5) node[above] {$\omega$};
  \node[below left,inner sep=2pt] at (0,0) {$0$};
  \foreach \c in {-1.9,1.9}{
    \filldraw[figoutline,fill=white] (\c,0) ellipse[x radius=.50,y radius=.34];
  }
  \node at (-1.9,0) {$-K$};
  \node at (1.9,0) {$K$};
  \node[above=4pt] at (-1.9,1.3) {$\Delta_-\times\R^d$};
  \node[above=4pt] at (1.9,1.3) {$\Delta_+\times\R^d$};
  \fill (2.30,.82) circle[radius=1.15pt];
  \node[right=3pt] at (2.30,.82) {$z_n$};
  \draw[figleader] (-1.9,-.38)--(-1.9,-.49);
  \draw[figleader] (1.9,-.38)--(1.9,-.49);
  \node[below=2pt] at (-1.9,-.49) {$-s_0$};
  \node[below=2pt] at (1.9,-.49) {$s_0$};
\end{tikzpicture}
\caption{The prescribed holes $K$ and $-K$ inside the two delay strips.
Points $z_n$ are chosen in $(\Delta_+\times\R^d)\setminus K$.
Hermitian symmetry gives the corresponding support inclusion in the negative strip.}
\label{fig:inner-support}
\end{figure}
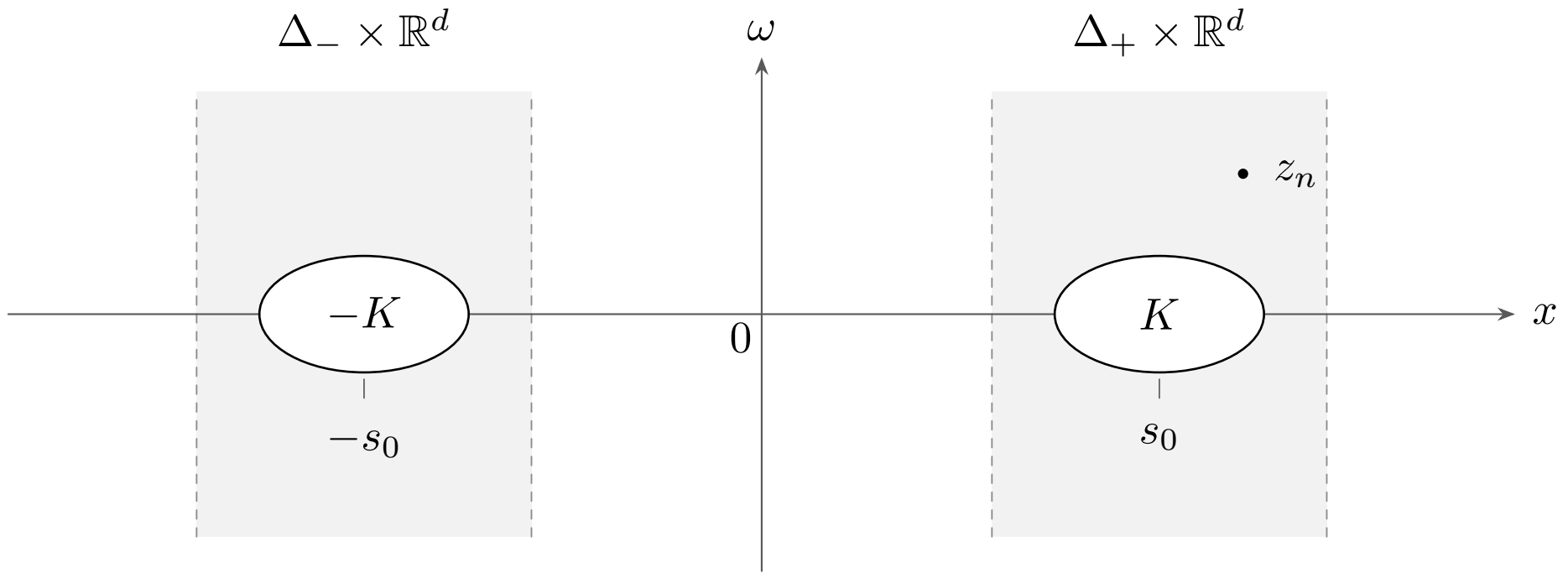

\subsection{Admissible translations}
\label{subsec:corridor}

The function $H$ defined in \eqref{eq:H} satisfies $A(H)=0$ on $K\cup-K$ as well as the inclusions \eqref{eq:inner-support}.  We next
choose translations $\alpha$ such that cubes based at $\alpha$ do not
contribute for $x\in\Delta$.  We also require both $\alpha$ and $\alpha+q$ to be
available for any prescribed separation $q$.

The following measure theoretic argument
provides a periodic set of translations with these properties. Here, we let $\pi:\R^d\to\R^d/\Z^d$ be the quotient map, and equip the torus $\R^d/\Z^d$
with Haar measure of total mass $1$.

\begin{lemma}
\label{lem:corridor}
Let $U,E\subset\R^d$ be compact and suppose that $\vol(U+E)<1/4$.
Then there is an open set $R_0\subset\R^d/\Z^d$ with Haar measure
$|R_0|>1/2$ such that, with
$\mathcal P=\pi^{-1}(R_0)$, one has
\begin{equation}
 (\mathcal P-(\Lambda+U))\cap E=\varnothing.
 \label{eq:corridor-safety-abstract}
\end{equation}
Moreover, for every $y\in\R^d$, the set
\[
 \mathcal A_y=\{\alpha\in\mathcal P:\alpha+y\in\mathcal P\}
\]
contains infinitely many pairwise disjoint translates of a nonempty open
ball.  Consequently, deleting a finite union of bounded sets from
$\mathcal A_y$ leaves a nonempty set.
\end{lemma}

\begin{proof}
Define the compact set $B=\pi(U+E)\cup\pi(\vartheta+U+E)\subset\R^d/\Z^d$ and let $R_0=(\R^d/\Z^d)\setminus B$. Then $R_0$ is open, and since $|B|\le2\vol(U+E)<1/2$, we have $|R_0|>1/2$.

To prove \eqref{eq:corridor-safety-abstract}, suppose by contradiction that there exist $\alpha\in\mathcal P$, $\lambda\in\Lambda$ and $u\in U$ such that $\alpha-(\lambda+u)\in E$. Then $\alpha\in\lambda+U+E$, and since $\lambda\in\Z^d\cup(\Z^d+\vartheta)$, it follows that $\pi(\alpha)\in B$, which contradicts $\pi(\alpha)\in R_0$.

Now let $y\in\R^d$. The translation invariance of the Haar measure gives
$$
|R_0\cap(R_0-\pi(y))|\ge2|R_0|-1>0.
$$
Hence, $\mathcal A_y=\pi^{-1}\bigl(R_0\cap(R_0-\pi(y))\bigr)$ is a nonempty open set satisfying $\mathcal A_y+\Z^d=\mathcal A_y$. In particular, $\mathcal A_y$ contains an open ball $G$ of radius less than $1/2$ together with the pairwise disjoint translates $G+n$, $n\in\Z^d$. Since every bounded set intersects only finitely many of these translates, the last claim follows.
\end{proof}

Set
\[
 E_H:=(s_0e_1+[-(L+2\rho),L+2\rho]^d)\cup (-s_0e_1+[-(L+2\rho),L+2\rho]^d),
 \]
and
\[
 U:=\{0\}\cup (s_0e_1+[-(L+3\rho),L+2\rho]^d).
\]
Our fixed parameters from \eqref{eq:parameters} give $\mathrm{vol}(U+E_H) \leq 2(2L+4\rho)^d + 2(4L+9\rho)^d < 1/4$, so we may apply Lemma~\ref{lem:corridor}.
Let
$\mathcal P$ be the resulting periodic set.

By the definition of $W$, $V^\Psi$ and $R_n$ (see Lemma~\ref{lem:local-packet}), the function $H$ is the limit in $L^2(\R^d)$ of finite linear combinations of the functions $\tau_s\1_{L\qq}$ with $s\in\Lambda+U$. This fact is used in the proof of the following lemma.

\begin{lemma}
\label{lem:corridor-H}
If $C=\alpha+[0,\varepsilon]^d$, where
$\alpha\in\mathcal P$ and $0<\varepsilon\le\rho$, then
\[
 A(\1_C,H)(x,\omega)=A(H,\1_C)(x,\omega)=0
 \qquad\text{for }x\in \Delta.
\]
\end{lemma}

\begin{proof}
Approximate $H$ in $L^2$ by finite linear combinations of the functions
$\tau_s\1_{L\qq}$ with $s\in\Lambda+U$.  If
$A(\1_C,\tau_s\1_{L\qq})(x,\omega)$ were nonzero for some $x\in\Delta$, then
\eqref{eq:support-rule} would give
\[
 x\in(\alpha-s)+[-L,\varepsilon]^d,
\]
and hence $\alpha-s\in E_H$.  This contradicts equation \eqref{eq:corridor-safety-abstract}.

On the other hand, nonvanishing would give $s-\alpha\in E_H$ and since
$E_H=-E_H$ we arrive at the same contradiction. The conclusion for $H$ follows
after passing to the $L^2$ limit.
\end{proof}

\subsection{Support at the remaining delays}
\label{subsec:outer-support}

A pair of indicators of small cubes whose starting points differ by $y$ has
a mixed ambiguity term centered at delay $y$.  For each point in a countable
dense set, we place such a pair in $\mathcal P$, far enough from the earlier
pairs that no new term contributes on $\Delta$.  Lemma~\ref{lem:diagonal}
will give nonvanishing on that dense set, and hence the required support
inclusion.

Choose a countable dense set $\mathcal D\subset\R^d\times\R^d$ of the form
\[
 \mathcal D=\{\zeta_m=(q_m,\omega_m):m\ge1\}
 \subset\bigl((\R^d\setminus\overline\Delta)\setminus\{0\}\bigr)
       \times\R^d.
\]
Since $\Delta$ is regular closed, the closure of $\mathcal D$ satisfies
\[
 \overline{\mathcal D}=(\R^d\setminus\mathrm{int}\Delta)\times\R^d.
\]
For each $\zeta_m$ we will construct functions
\[
 \psi_m
 =\1_{\alpha_m+[0,\varepsilon_m]^d}
  +\1_{\beta_m+[0,\varepsilon_m]^d},
 \qquad \beta_m-\alpha_m=q_m,
\]
with $\alpha_m,\beta_m\in\mathcal P$ and
$0<\varepsilon_m\le\rho$.

We construct the functions $\psi_m$ by induction on $m$. Suppose that $\psi_1,\ldots,\psi_{m-1}$ have been constructed, and define the sets $S_{m-1}=\{\alpha_j,\beta_j:1\le j\le m-1\}$ and $E_g=\Delta+[-\rho,\rho]^d$. By Lemma~\ref{lem:corridor}, there exists
$$
\alpha_m\in\mathcal A_{q_m}\setminus\bigcup_{\gamma\in S_{m-1}}\bigl((\gamma+E_g)\cup(\gamma-q_m+E_g)\bigr),
$$
and we set $\beta_m=\alpha_m+q_m$. Since $\Delta$ is closed with $\Delta=-\Delta$ and $0,q_m\notin\Delta$, we can choose $0<\varepsilon_m\le\rho$ with $\varepsilon_m<\|q_m\|_\infty$ such that
\begin{equation}\label{eq:eps-choice}
\bigl([-\varepsilon_m,\varepsilon_m]^d\cup(q_m+[-\varepsilon_m,\varepsilon_m]^d)\cup(-q_m+[-\varepsilon_m,\varepsilon_m]^d)\bigr)\cap\Delta=\varnothing
\end{equation}
and such that $\sin(\pi(\omega_m)_j\varepsilon_m)\neq0$ for every $j\in\{1,\dots,d\}$ with $(\omega_m)_j\neq0$.

We proceed by computing $A(\psi_m)$. Writing $G_m=A(\1_{\varepsilon_m\qq})$, Lemma~\ref{lem:basic-facts}\ref{fact:translation} gives
$$
\begin{aligned}
A(\psi_m)(x,\omega)&=\bigl(e^{-2\pi i\omega\cdot\alpha_m}+e^{-2\pi i\omega\cdot\beta_m}\bigr)G_m(x,\omega)\\
&\quad+e^{-\pi i\omega\cdot(\alpha_m+\beta_m)}\bigl(G_m(x-q_m,\omega)+G_m(x+q_m,\omega)\bigr),
\end{aligned}
$$
and by Lemma~\ref{lem:basic-facts}\ref{fact:expra1q}, we have $G_m(x,\omega)=0$ for $x\notin[-\varepsilon_m,\varepsilon_m]^d$. Hence, \eqref{eq:eps-choice} shows that $A(\psi_m)=0$ on $\Delta\times\R^d$. Moreover, since $\varepsilon_m<\|q_m\|_\infty$, we have $G_m(q_m,\omega_m)=G_m(2q_m,\omega_m)=0$, and therefore
$$
A(\psi_m)(\zeta_m)=e^{-\pi i\omega_m\cdot(\alpha_m+\beta_m)}G_m(0,\omega_m)=e^{-\pi i\omega_m\cdot(\alpha_m+\beta_m+\varepsilon_m\one)}\prod_{j=1}^d\varepsilon_m\sinc\bigl(\pi(\omega_m)_j\varepsilon_m\bigr)\neq0.
$$

Next, Lemma~\ref{lem:corridor-H} gives $A(\psi_m,H)=A(H,\psi_m)=0$ on $\Delta\times\R^d$. Further, let $j<m$, $\gamma\in\{\alpha_m,\beta_m\}$ and $\gamma'\in\{\alpha_j,\beta_j\}$. By the choice of $\alpha_m$, we have $\gamma-\gamma'\notin E_g$, and since $E_g=-E_g$ and $\varepsilon_j,\varepsilon_m\le\rho$, the sets $(\gamma-\gamma')+[-\varepsilon_j,\varepsilon_m]^d$ and $(\gamma'-\gamma)+[-\varepsilon_m,\varepsilon_j]^d$ are disjoint from $\Delta$. By \eqref{eq:support-rule}, it follows that $A(\psi_m,\psi_j)=A(\psi_j,\psi_m)=0$ on $\Delta\times\R^d$. Consequently, for all $N\in\N$ and $s_1,\dots,s_N\in\R$ we have
\begin{equation}\label{eq:finite-sum}
A\Bigl(H+\sum_{j=1}^Ns_j\psi_j\Bigr)=A(H)\quad\text{on }\Delta\times\R^d.
\end{equation}

We now apply Lemma~\ref{lem:diagonal} with $h=H$, $g_m=\psi_m$, $z_m=\zeta_m$ and $E=K\cup-K$. Since $A(H)=0$ on $K\cup-K\subset\Delta\times\R^d$, condition~\ref{hyp:diagonal-1} of Lemma~\ref{lem:diagonal} follows from \eqref{eq:finite-sum}. Condition~\ref{hyp:diagonal-2} is satisfied, since for every $g\in L^2(\R^d)$ the coefficient of $s^2$ in the polynomial $s\mapsto A(H+g+s\psi_m)(\zeta_m)$ is given by $A(\psi_m)(\zeta_m)\neq0$. Hence, there exist real numbers $s_m$ such that the series
\begin{equation}\label{eq:f-def}
f=H+\sum_{m=1}^\infty s_m\psi_m
\end{equation}
converges in $L^2(\R^d)$ and satisfies
$$
A(f)=0\quad\text{on }K\cup-K,\qquad A(f)(\zeta_m)\neq0,\quad m\ge1.
$$
Moreover, \eqref{eq:finite-sum} and Lemma~\ref{lem:basic-facts}\ref{fact:continuity} give
\begin{equation}\label{eq:f-equals-H}
A(f)=A(H)\quad\text{on }\Delta\times\R^d.
\end{equation}

Since every $\zeta_m$ belongs to the closed set $\supp A(f)$ and $\overline{\mathcal D}=(\R^d\setminus\interior\Delta)\times\R^d$, it follows that
\begin{equation}\label{eq:outer-support}
(\R^d\setminus\Delta)\times\R^d\subset\supp A(f).
\end{equation}

\begin{proof}[Proof of the local realization]
Let $f$ be the function defined in \eqref{eq:f-def}. By \eqref{eq:f-equals-H}, we have $A(f)(z_n)=A(H)(z_n)\neq0$ for every $n\ge1$. Since $Z_+=\{z_n:n\ge1\}$ is dense in $(\Delta_+\times\R^d)\setminus K$ and $\supp A(f)$ is a closed set satisfying $\supp A(f)=-\supp A(f)$, it follows that
$$
(\Delta_+\times\R^d)\setminus K\subset\supp A(f),\qquad(\Delta_-\times\R^d)\setminus-K\subset\supp A(f).
$$
Together with \eqref{eq:outer-support}, this gives $\R^{2d}\setminus(K\cup-K)\subset\supp A(f)$.

Since $K$ and $-K$ are disjoint compact sets and $\Omega$ is regular open, the interior of $K\cup-K$ is given by $\Omega\cup-\Omega$, and since $\supp A(f)$ is closed, we obtain
$$
\R^{2d}\setminus(\Omega\cup-\Omega)=\overline{\R^{2d}\setminus(K\cup-K)}\subset\supp A(f).
$$
Conversely, $A(f)$ vanishes on the open set $\Omega\cup-\Omega$, and therefore $\supp A(f)\subset\R^{2d}\setminus(\Omega\cup-\Omega)$. This shows that $\supp A(f)=\R^{2d}\setminus(\Omega\cup-\Omega)$, that is, $\Zopen(f)=\Omega\cup-\Omega$.
\end{proof}

\section{Annular obstructions}

\subsection{A Fourier obstruction}
\label{sec:obstruction}

By Theorem~\ref{thm:main}, sufficiently small rescalings of bounded regular
open sets, placed away from the origin and symmetrized, occur as open zero
sets. Hence, obstructions to open zero sets can only appear at larger scales,
and they require a global property of ambiguity functions. The property we use
is the invariance of $|A(f)|^2$ under the symplectic Fourier transform. It
implies that $|A(f)|^2$ is orthogonal to every real function $h$ satisfying
$\mathcal F_\sigma h=-h$, which prevents the negative set of such a function
$h$ from being contained in an open zero set
(Proposition~\ref{prop:certificate}).

\begin{lemma}
\label{lem:self}
For every $f\in L^2(\R^d)$, it holds that
\[
 F_f:=|A(f)|^2\in L^1(\R^{2d})\cap L^2(\R^{2d})
 \quad\text{and}\quad
 \mathcal F_\sigma F_f=F_f.
\]
\end{lemma}

\begin{proof}
Since the ambiguity function and the short-time Fourier transform are related by
\[
 V_ff(x,\omega)=e^{-\pi i x\cdot\omega}A(f)(x,\omega),
\]
we have $|V_ff|=|A(f)|$.  Applying
Lemma~\ref{lem:spectrogram-identity} with $w=f$ gives
\[
 \mathcal F_\sigma(|A(f)|^2)=A(f)\overline{A(f)}=|A(f)|^2.
\]
Finally, since $A(f)\in L^2(\R^{2d})\cap\mathcal{C}_0(\R^{2d})$, we have $|A(f)|^2\in L^1(\R^{2d})\cap L^\infty(\R^{2d})\subset L^2(\R^{2d})$.
\end{proof}

The identity $\mathcal F_\sigma(|A(f)|^2)=|A(f)|^2$ in Lemma~\ref{lem:self} is classical in the radar literature \cite{AuslanderTolimieri1985,Wilcox1991}. This identity gives a sign obstruction.  If $h$ is a real Schwartz function
and $\mathcal F_\sigma h=-h$, then $h$ is orthogonal to $|A(f)|^2$.

\begin{proposition}
\label{prop:certificate}
Let $h\in\Sclass(\R^{2d})$ be real valued and satisfy
$\mathcal F_\sigma h=-h$ and $h(0)>0$. Let $O\subset\R^{2d}$ be open, with
$0\notin\overline O$, and assume that $\{z:h(z)<0\}\subset O$. Then no
nonzero $f\in L^2(\R^d)$ can satisfy $O\subset\Zopen(f)$.
\end{proposition}

\begin{proof}
Suppose by contradiction that $O\subset\Zopen(f)$, and set
$F=|A(f)|^2$.  By Lemma~\ref{lem:self}, we have
$\mathcal F_\sigma F=F$, and an application of Plancherel's formula gives
\[
 \int_{\R^{2d}}F(z)h(z)\,dz
 =\int_{\R^{2d}}(\mathcal F_\sigma F)(z)
   \overline{(\mathcal F_\sigma h)(z)}\,dz
 =-\int_{\R^{2d}}F(z)h(z)\,dz.
\]
Thus
\begin{equation}
 \int_{\R^{2d}}Fh=0.
 \label{eq:orthogonality}
\end{equation}

On the set $\{z\in\R^{2d}:h(z)<0\}$, the function $F$ vanishes because
$\{h<0\}\subset O\subset\Zopen(f)$.  Hence $Fh\ge0$ everywhere.
At the origin, we have
\[
 F(0)=\|f\|_2^4>0
 \quad\text{and}\quad
 h(0)>0.
\]
By continuity, $Fh$ is strictly positive on a neighborhood of the origin,
contradicting \eqref{eq:orthogonality}.
\end{proof}

\begin{remark}
\label{rem:sign-uncertainty}
Proposition~\ref{prop:certificate} belongs to a familiar family of sign-based certificates in Fourier analysis.  Sign conditions on a function and its Fourier transform underlie the linear programming bounds for sphere packings of Cohn and Elkies \cite{CohnElkies2003} as well as the sign uncertainty principles of Bourgain, Clozel, and Kahane \cite{BourgainClozelKahane2010}; in the linear programming setting, the radial Laguerre--Gaussian eigenfunctions employed in Section~\ref{sec:annuli} serve as basic building blocks \cite[p.~704]{CohnElkies2003}.  

In our case, the additional property used is that $|A(f)|^2$ is a \emph{fixed point} of the symplectic
Fourier transform, so that every real Schwartz function $h$ with
$\mathcal F_\sigma h=-h$ and $h(0)>0$ produces a "forbidden" region
$\{h<0\}$ for open zero sets.
\end{remark}

\subsection{Laguerre--Gaussian functions}
\label{sec:annuli}

Radial Laguerre-Gaussian functions give explicit examples of certificates because their Fourier eigenvalue is determined by the parity of the Laguerre index.

To make use of this property, let $L_k^\alpha$ denote the generalized Laguerre polynomial, and define
\[
 \psi_k^{(d)}(z)
 :=L_k^{d-1}(2\pi|z|^2)e^{-\pi|z|^2},
 \qquad z\in\R^{2d}.
\]
On $\R^{2d}$, we use the Fourier transform convention
\[
 \widehat F(\xi)=\int_{\R^{2d}}F(z)e^{-2\pi i\xi\cdot z}\,dz.
\]

\begin{lemma}
\label{lem:laguerre}
For every integer $k\geq0$, one has
$\widehat{\psi_k^{(d)}}=(-1)^k\psi_k^{(d)}$ and
$\mathcal F_\sigma\psi_k^{(d)}=(-1)^k\psi_k^{(d)}$.
\end{lemma}

\begin{proof}
The functions $\psi_k^{(d)}$ form, up to normalization, a basis of the
radial eigenfunctions of the Fourier transform of $\R^{2d}$. With our
normalization, the eigenvalue formula
$\widehat{\psi_k^{(d)}}=(-1)^k\psi_k^{(d)}$ is recorded, for instance,
in \cite[p.~704]{CohnElkies2003} (there with $n=2d$ and $\alpha=d-1$),
and it goes back to classical work on Laguerre functions
\cite[\S4.23]{Lebedev1972}.  The symplectic statement follows from
$\mathcal F_\sigma F(x,\omega)=\widehat F(\omega,-x)$ and the fact that a radial
function takes the same value at $(\omega,-x)$ and $(x,\omega)$.
\end{proof}

\subsection{Dimension one}
\label{subsec:dimension-one}

We seek a combination of odd Laguerre functions whose radial polynomial is
positive at the origin and negative on one bounded interval.

In dimension
$d=1$, write $L_k=L_k^0$ and set
\begin{equation}
 h(z)
 =\Bigl(17L_1(2\pi|z|^2)-6L_3(2\pi|z|^2)\Bigr)
  e^{-\pi|z|^2},
 \qquad z\in\R^2.
 \label{eq:h-explicit}
\end{equation}
Both Laguerre indices are odd, so
$\mathcal F_\sigma h=-h$ by Lemma~\ref{lem:laguerre}. Information about the sign of $h$ can be obtained dicrectly from the elementary identities
\[
 L_1(x)=1-x,
 \qquad
 L_3(x)=1-3x+\frac32x^2-\frac16x^3,
\]
which give the relation
\begin{equation}
 17L_1(x)-6L_3(x)
 =(x+1)\bigl(x-(5-\sqrt{14})\bigr)
       \bigl(x-(5+\sqrt{14})\bigr).
 \label{eq:factorization}
\end{equation}

\begin{proof}[Proof of Theorem~\ref{thm:counterexample}]
Take $r=\sqrt{(5-\sqrt{14})/(2\pi)}$ and $R=\sqrt{(5+\sqrt{14})/(2\pi)}$ and set
$O_*=\{z\in\R^2\,: r<|z|<R\}$. Then $O_*$ is bounded, regular open, centrally
symmetric, and separated from the origin.
By \eqref{eq:factorization}, the negative set of the function $h$ in
\eqref{eq:h-explicit} is exactly the annulus $O_*$.  Moreover,
$h(0)=11>0$ and $\mathcal F_\sigma h=-h$.
Proposition~\ref{prop:certificate} therefore shows that $O_*$ cannot be
contained in $\Zopen(f)$ for any nonzero $f$.
\end{proof}

Symplectic changes of variables give nonradial examples as well.

\begin{corollary}
\label{cor:symplectic}
For every $S\in\Sp(2,\R)=\mathrm{SL}(2,\R)$, the elliptic annulus $S(O_*)$
cannot be contained in $\Zopen(f)$ for any nonzero $f\in L^2(\R)$.
\end{corollary}

\subsection{Higher dimensions}
\label{subsec:higher-dimensions}

We choose the coefficients in the following definition of $q_d$ such that $q_d$ is a monic polynomial of degree $3$ satisfying $q_d(0)>0$ and $q_d(d+1)<0$.  Define
\[
 h_d(z)=q_d(2\pi|z|^2)e^{-\pi|z|^2},
\]
where
\begin{align*}
 q_d(x)
 &:=\bigl((d+1)(d+2)+1\bigr)L_1^{d-1}(x)-6L_3^{d-1}(x)\\
 &=x^3-3(d+2)x^2
   +\bigl(2(d+1)(d+2)-1\bigr)x+d.
\end{align*}
Again both Laguerre indices are odd, so
$\mathcal F_\sigma h_d=-h_d$.  The polynomial $q_d$ satisfies
\[
 q_d(0)=d>0,
 \qquad
 q_d(d+1)=-(d^2+2d+2)<0,
 \qquad
 q_d(x)\longrightarrow+\infty
 \quad(x\to+\infty).
\]
The intermediate value theorem gives a root $a_d\in(0,d+1)$ and a root
$b_d\in(d+1,\infty)$.  These roots are distinct.  Denoting by $c_d$ the third root of $q_d$, Vieta's formulas give $a_db_dc_d=-d<0$, and hence $c_d<0$. Since $q_d(x)=(x-c_d)(x-a_d)(x-b_d)$, it follows that for $x\ge0$ we have the equivalence
$$
q_d(x)<0\quad\Longleftrightarrow\quad a_d<x<b_d.
$$

\begin{corollary}
\label{cor:all-dimensional-annuli}
For every $d\ge1$, the set
\[
 O_d=\{z\in\R^{2d}:a_d<2\pi|z|^2<b_d\}
\]
is bounded, regular open, centrally symmetric, does not contains the origin, and has the following property: there does not exist non-zero $f\in L^2(\R^d)$ such that $O_d \subseteq \Zopen(f)$. The same conclusion holds for every symplectic image
$S(O_d) = \{ So : o \in O_d \}$ with $S\in\Sp(2d,\R)$.
\end{corollary}

\begin{proof}
The negative set of $h_d$ is exactly $O_d$, and we have
$h_d(0)=d>0$ and $\mathcal F_\sigma h_d=-h_d$.
Proposition~\ref{prop:certificate} gives the first assertion.  For the
second, apply the same proposition to $h_d\circ S^{-1}$.
\end{proof}

\section{Phase retrieval}
\label{sec:phase-retrieval}

The key point in the proof of Theorem~\ref{thm:phase-retrieval-d} is that
two rank-one ambiguity functions cannot differ only on a compact subset of
phase space unless the underlying functions are linearly dependent.

\begin{lemma}
\label{lem:compact-rigidity}
Let $f,g\in L^2(\R^d)$.  If $A(f)-A(g)$ has compact support in
$\R^{2d}$, then $f$ and $g$ are linearly dependent.
\end{lemma}

\begin{proof}
We start by defining $D=A(f)-A(g)$ and $K(s,t)=f(s)\overline{f(t)}-g(s)\overline{g(t)}$.
Let $C\subset\R^d$ be a compact set containing the projection of
$\supp D$ onto the delay variable.  If $x\notin C$, then
\[
 0=D(x,\omega)
 =\int_{\R^d}
 K\!\left(t+\frac{x}{2},t-\frac{x}{2}\right)
 e^{-2\pi i\omega\cdot t}\,dt,
 \qquad\omega\in\R^d.
\]
The integrand belongs to $L^1(\R^d)$.  Uniqueness of the Fourier transform,
followed by Fubini, therefore gives
\begin{equation}
 K(s,t)=0
 \quad\text{for almost every }(s,t)\text{ satisfying }s-t\notin C.
 \label{eq:kernel-band-d}
\end{equation}

Now suppose that $f$ and $g$ are linearly independent.  For a compactly
supported $\varphi\in L^2(\R^d)$, define
\[
 G_\varphi(s)=\int_{\R^d}K(s,t)\varphi(t)\,dt.
\]
By \eqref{eq:kernel-band-d}, the function $G_\varphi$ is supported in
the compact set $\supp\varphi+C$.  On the other hand,
\[
 G_\varphi(s)
 =f(s)\int_{\R^d}\overline{f(t)}\varphi(t)\,dt
  -g(s)\int_{\R^d}\overline{g(t)}\varphi(t)\,dt.
\]
The linear map
$$
\varphi\longmapsto\Bigl(\int_{\R^d}\overline{f(t)}\varphi(t)\,dt,\ \int_{\R^d}\overline{g(t)}\varphi(t)\,dt\Bigr)
$$
from $L^2(\R^d)$ to $\C^2$ is continuous, and it is surjective since $f$ and $g$ are linearly independent. Since the compactly supported functions are dense in $L^2(\R^d)$, their image under this map is a dense subspace of $\C^2$ and hence equal to $\C^2$. Thus, there exist compactly supported $\varphi_1,\varphi_2\in L^2(\R^d)$ such that
\[
 G_{\varphi_1}=f,
 \qquad
 G_{\varphi_2}=-g.
\]
It follows that $f$ and $g$ are compactly supported, and therefore $\widehat f$ and $\widehat g$ extend to entire functions on $\C^d$.

By \eqref{eq:cov-fourier}, the function $A(\widehat f)-A(\widehat g)$ has compact support, and applying the above argument to the linearly independent functions $\widehat f$ and $\widehat g$ shows that $\widehat f$ and $\widehat g$ are compactly supported. Hence, $\widehat f=\widehat g=0$, which contradicts the linear independence of $f$ and $g$.
\end{proof}

\begin{remark}
\label{rem:support-uncertainty}
The special case $g=0$ of Lemma~\ref{lem:compact-rigidity} shows that a nonzero ambiguity function is never compactly supported. 
As mentioned in the introduction, this result is well known and can even be extended
to sets of finite measure (and many more).
The two-function statement of Lemma~\ref{lem:compact-rigidity}
also follows from the finite-rank Weyl-transform theorem
\cite[Theorem~2.2]{NarayananRatnakumar2010}, applied to $A(f)-A(g)$,
whose corresponding operator has rank at most two.
We include the elementary proof above for completeness.
\end{remark}

\begin{proof}[Proof of Theorem~\ref{thm:phase-retrieval-d}]
Suppose that $|V_wf|=|V_wg|$, and set $D=A(f)-A(g)$.
Lemma~\ref{lem:spectrogram-identity} gives
\[
 D\,\overline{A(w)}=0.
\]
At every point of $\supp A(w)$, each neighborhood contains a point at which
$A(w)$ is nonzero.  The preceding identity shows that $D$ vanishes at all
such points, and continuity then gives
\[
 D=0\qquad\text{on }\supp A(w).
\]
Consequently,
\begin{equation}
 \supp D\subset\overline{\Zopen(w)}.
 \label{eq:supp-D}
\end{equation}
The right-hand side of \eqref{eq:supp-D} is compact by hypothesis, so
Lemma~\ref{lem:compact-rigidity} implies that $f$ and $g$ are linearly
dependent.

Since $A(w)(0)=\|w\|_2^2>0$, the identity
$D\,\overline{A(w)}=0$ also gives
$\|f\|_2^2=A(f)(0)=A(g)(0)=\|g\|_2^2$.
If $f=0$, the equality of norms gives $g=0$.  Otherwise $g=cf$ for some
$c\in\C$, and the same equality gives $|c|=1$.  Hence
$g=e^{i\theta}f$ for some $\theta\in\R$.
\end{proof}

Every open zero set produced by Theorem~\ref{thm:main} is bounded.  Thus each
window constructed there satisfies the injectivity conclusion of
Theorem~\ref{thm:phase-retrieval-d}.

\begin{remark}
\label{rem:finite-measure}
Theorem~\ref{thm:phase-retrieval-d} remains true if $\Zopen(w)$ is only assumed to have finite measure. Indeed, if $|V_wf|=|V_wg|$, then $D=A(f)-A(g)$ vanishes outside $\Zopen(w)$ by the proof of Theorem~\ref{thm:phase-retrieval-d}, and it suffices to show that $D=0$. By Moyal's identity, $\mathcal M=\{A(f,u)+A(g,v):u,v\in L^2(\R^d)\}$ is a closed subspace of $L^2(\R^{2d})$ with reproducing kernel $k_z$ satisfying $\|k_z\|_2^2\le2$. Hence, as in \cite{AmreinBerthier,Benedicks1985}, the functions in $\mathcal M$ that vanish outside a set $S$ of finite measure form a space of dimension at most $2|S|$. If $D\neq0$, then the functions
$$
D_a(x,\omega)=e^{\pi i\omega\cdot a}D(x-a,\omega)=\bigl(A(f,\tau_{-a}f)-A(g,\tau_{-a}g)\bigr)(x,\omega),\quad a\in\R^d,
$$
belong to $\mathcal M$, vanish outside $\Zopen(w)+(a,0)$, and are linearly independent for distinct $a$ (since distinct translates of a nonzero function in $L^2(\R^d)$ are linearly independent). For distinct and sufficiently small $a_1,\dots,a_N\in\R^d$, the set $S=\bigcup_{j=1}^N(\Zopen(w)+(a_j,0))$ satisfies $|S|\le|\Zopen(w)|+1$, and therefore $N\le2(|\Zopen(w)|+1)$ for every $N\in\N$, which gives a contradiction.
\end{remark}

\section{Appendix: Lean formalization}
\label{sec:lean}

This appendix describes the Lean formalization of the three main results of
the paper, Theorems~\ref{thm:phase-retrieval-d}, \ref{thm:main}
and~\ref{thm:counterexample}. The formalization is written in Lean~4
\cite{lean4} and builds on Mathlib \cite{mathlib}, the mathematical library of
Lean. Its source code is available at
\begin{center}
\url{https://github.com/lukasliehr/Ambiguity-Zeros}.
\end{center}
The file \mbox{\ttfamily\frenchspacing Showcase\textcolor[rgb]{0.20,0.40,0.64}{.}lean} therein provides a self-contained version of the
three statements that were verified, together with all definitions that are
not already available in Mathlib. Its only import is Mathlib. The statements
in \mbox{\ttfamily\frenchspacing Showcase\textcolor[rgb]{0.20,0.40,0.64}{.}lean} have been carefully curated and reviewed by the authors
with the objective of making them understandable to a broad audience. The
remaining Lean code, which contains the proofs, has been generated with the
assistance of large language models (mainly Codex~5.4 and Codex~5.5).

The formalization relies on many mathematical objects that are already
defined in Mathlib, such as the real and complex numbers, the Lebesgue measure,
the Bochner integral, the spaces $L^p$, and the interior and the closure of a
set. These objects have been carefully designed and reviewed by the Mathlib
community, and we do not discuss their construction. The objective of this
appendix is to make the correspondence between the Lean statements and
Theorems~\ref{thm:phase-retrieval-d}--\ref{thm:counterexample} precise and
transparent without assuming familiarity with Lean.

\subsection{Euclidean space, \texorpdfstring{$L^2$}{L2} and phase space}

All declarations of this appendix are placed inside the namespace
\mbox{\ttfamily\frenchspacing AmbiguityHoles}; for readability, the namespace prefix is omitted in what
follows. The underlying spaces are introduced as abbreviations.

\par\medskip
{\ttfamily\frenchspacing\setlength{\parindent}{0pt}\setlength{\parskip}{0pt}\raggedright
\noindent \textcolor[rgb]{0.13,0.29,0.53}{\textbf{abbrev}} Rd \textcolor[rgb]{0.81,0.36,0.00}{\textbf{(}}d \textcolor[rgb]{0.81,0.36,0.00}{\textbf{:}} \makebox[0.525em][c]{$\mathbb{N}$}\textcolor[rgb]{0.81,0.36,0.00}{\textbf{)}} \textcolor[rgb]{0.81,0.36,0.00}{\textbf{:=}} EuclideanSpace \makebox[0.525em][c]{$\mathbb{R}$} \textcolor[rgb]{0.81,0.36,0.00}{\textbf{(}}Fin d\textcolor[rgb]{0.81,0.36,0.00}{\textbf{)}}\par
\noindent \textcolor[rgb]{0.13,0.29,0.53}{\textbf{abbrev}} L2 \textcolor[rgb]{0.81,0.36,0.00}{\textbf{(}}d \textcolor[rgb]{0.81,0.36,0.00}{\textbf{:}} \makebox[0.525em][c]{$\mathbb{N}$}\textcolor[rgb]{0.81,0.36,0.00}{\textbf{)}} \textcolor[rgb]{0.81,0.36,0.00}{\textbf{:=}} Lp \makebox[0.525em][c]{$\mathbb{C}$} \textcolor[rgb]{0,0,0.81}{\textbf{2}} \textcolor[rgb]{0.81,0.36,0.00}{\textbf{(}}volume \textcolor[rgb]{0.81,0.36,0.00}{\textbf{:}} Measure \textcolor[rgb]{0.81,0.36,0.00}{\textbf{(}}Rd d\textcolor[rgb]{0.81,0.36,0.00}{\textbf{)}}\textcolor[rgb]{0.81,0.36,0.00}{\textbf{)}}\par
\noindent \textcolor[rgb]{0.13,0.29,0.53}{\textbf{abbrev}} PhaseSpace \textcolor[rgb]{0.81,0.36,0.00}{\textbf{(}}d \textcolor[rgb]{0.81,0.36,0.00}{\textbf{:}} \makebox[0.525em][c]{$\mathbb{N}$}\textcolor[rgb]{0.81,0.36,0.00}{\textbf{)}} \textcolor[rgb]{0.81,0.36,0.00}{\textbf{:=}} Rd d \textcolor[rgb]{0.20,0.40,0.64}{\makebox[0.525em][c]{$\times$}} Rd d\par
}\par\medskip

The type \mbox{\ttfamily\frenchspacing Fin d} consists of the indices $0,\ldots,d-1$, and
\mbox{\ttfamily\frenchspacing EuclideanSpace \makebox[0.525em][c]{$\mathbb{R}$} \textcolor[rgb]{0.81,0.36,0.00}{\textbf{(}}Fin d\textcolor[rgb]{0.81,0.36,0.00}{\textbf{)}}} is the space $\R^d$ equipped with the
Euclidean inner product and norm. In Mathlib, \mbox{\ttfamily\frenchspacing volume} denotes the
Lebesgue measure, and \mbox{\ttfamily\frenchspacing Lp \makebox[0.525em][c]{$\mathbb{C}$} \textcolor[rgb]{0,0,0.81}{\textbf{2}} volume} is the Hilbert space
$L^2(\R^d)=L^2(\R^d,\C)$. Precisely, the elements of \mbox{\ttfamily\frenchspacing L2 d} are
equivalence classes of square-integrable functions modulo equality almost
everywhere. Consequently, for \mbox{\ttfamily\frenchspacing w \textcolor[rgb]{0.81,0.36,0.00}{\textbf{:}} L2 d} the condition \mbox{\ttfamily\frenchspacing w \textcolor[rgb]{0.20,0.40,0.64}{\makebox[0.525em][c]{$\neq$}} \textcolor[rgb]{0,0,0.81}{\textbf{0}}}
means that $w$ does not vanish almost everywhere, and the identity
\mbox{\ttfamily\frenchspacing f \textcolor[rgb]{0.20,0.40,0.64}{=} g} of two elements of \mbox{\ttfamily\frenchspacing L2 d} means that $f=g$ almost
everywhere, in accordance with the rest of the paper. Finally,
\mbox{\ttfamily\frenchspacing PhaseSpace d} is the product $\R^d\times\R^d$. For a point
\mbox{\ttfamily\frenchspacing z \textcolor[rgb]{0.81,0.36,0.00}{\textbf{:}} PhaseSpace d}, the components \mbox{\ttfamily\frenchspacing z\textcolor[rgb]{0.20,0.40,0.64}{.}\textcolor[rgb]{0,0,0.81}{\textbf{1}}} and \mbox{\ttfamily\frenchspacing z\textcolor[rgb]{0.20,0.40,0.64}{.}\textcolor[rgb]{0,0,0.81}{\textbf{2}}} correspond
to the delay variable $x$ and the frequency variable $\omega$ of
$z=(x,\omega)$.

\subsection{Ambiguity function and short-time Fourier transform}

The ambiguity function is formalized directly from its definition in
Section~\ref{sec:preliminaries}.

\par\medskip
{\ttfamily\frenchspacing\setlength{\parindent}{0pt}\setlength{\parskip}{0pt}\raggedright
\noindent \textcolor[rgb]{0.13,0.29,0.53}{\textbf{def}} ambiguity \textcolor[rgb]{0.81,0.36,0.00}{\textbf{\{}}d \textcolor[rgb]{0.81,0.36,0.00}{\textbf{:}} \makebox[0.525em][c]{$\mathbb{N}$}\textcolor[rgb]{0.81,0.36,0.00}{\textbf{\}}} \textcolor[rgb]{0.81,0.36,0.00}{\textbf{(}}f \textcolor[rgb]{0.81,0.36,0.00}{\textbf{:}} L2 d\textcolor[rgb]{0.81,0.36,0.00}{\textbf{)}} \textcolor[rgb]{0.81,0.36,0.00}{\textbf{(}}z \textcolor[rgb]{0.81,0.36,0.00}{\textbf{:}} PhaseSpace d\textcolor[rgb]{0.81,0.36,0.00}{\textbf{)}} \textcolor[rgb]{0.81,0.36,0.00}{\textbf{:}} \makebox[0.525em][c]{$\mathbb{C}$} \textcolor[rgb]{0.81,0.36,0.00}{\textbf{:=}}\par
\noindent \hspace*{1.05em}\textcolor[rgb]{0.20,0.40,0.64}{\makebox[0.525em][c]{$\int$}} t\textcolor[rgb]{0.81,0.36,0.00}{\textbf{,}} f \textcolor[rgb]{0.81,0.36,0.00}{\textbf{(}}t \textcolor[rgb]{0.20,0.40,0.64}{+} \textcolor[rgb]{0.81,0.36,0.00}{\textbf{(}}\textcolor[rgb]{0,0,0.81}{\textbf{1}} \textcolor[rgb]{0.20,0.40,0.64}{/} \textcolor[rgb]{0,0,0.81}{\textbf{2}} \textcolor[rgb]{0.81,0.36,0.00}{\textbf{:}} \makebox[0.525em][c]{$\mathbb{R}$}\textcolor[rgb]{0.81,0.36,0.00}{\textbf{)}} \textcolor[rgb]{0.20,0.40,0.64}{\makebox[0.525em][c]{$\bullet$}} z\textcolor[rgb]{0.20,0.40,0.64}{.}\textcolor[rgb]{0,0,0.81}{\textbf{1}}\textcolor[rgb]{0.81,0.36,0.00}{\textbf{)}} \textcolor[rgb]{0.20,0.40,0.64}{*} star \textcolor[rgb]{0.81,0.36,0.00}{\textbf{(}}f \textcolor[rgb]{0.81,0.36,0.00}{\textbf{(}}t \textcolor[rgb]{0.20,0.40,0.64}{-} \textcolor[rgb]{0.81,0.36,0.00}{\textbf{(}}\textcolor[rgb]{0,0,0.81}{\textbf{1}} \textcolor[rgb]{0.20,0.40,0.64}{/} \textcolor[rgb]{0,0,0.81}{\textbf{2}} \textcolor[rgb]{0.81,0.36,0.00}{\textbf{:}} \makebox[0.525em][c]{$\mathbb{R}$}\textcolor[rgb]{0.81,0.36,0.00}{\textbf{)}} \textcolor[rgb]{0.20,0.40,0.64}{\makebox[0.525em][c]{$\bullet$}} z\textcolor[rgb]{0.20,0.40,0.64}{.}\textcolor[rgb]{0,0,0.81}{\textbf{1}}\textcolor[rgb]{0.81,0.36,0.00}{\textbf{)}}\textcolor[rgb]{0.81,0.36,0.00}{\textbf{)}} \textcolor[rgb]{0.20,0.40,0.64}{*}\par
\noindent \hspace*{2.1em}Complex\textcolor[rgb]{0.20,0.40,0.64}{.}exp \textcolor[rgb]{0.81,0.36,0.00}{\textbf{(}}\textcolor[rgb]{0.20,0.40,0.64}{-}\textcolor[rgb]{0,0,0.81}{\textbf{2}} \textcolor[rgb]{0.20,0.40,0.64}{*} Real\textcolor[rgb]{0.20,0.40,0.64}{.}pi \textcolor[rgb]{0.20,0.40,0.64}{*} Complex\textcolor[rgb]{0.20,0.40,0.64}{.}I \textcolor[rgb]{0.20,0.40,0.64}{*} \textcolor[rgb]{0.81,0.36,0.00}{\textbf{(}}inner \makebox[0.525em][c]{$\mathbb{R}$} z\textcolor[rgb]{0.20,0.40,0.64}{.}\textcolor[rgb]{0,0,0.81}{\textbf{2}} t \textcolor[rgb]{0.81,0.36,0.00}{\textbf{:}} \makebox[0.525em][c]{$\mathbb{R}$}\textcolor[rgb]{0.81,0.36,0.00}{\textbf{)}}\textcolor[rgb]{0.81,0.36,0.00}{\textbf{)}}\par
}\par\medskip

Here \mbox{\ttfamily\frenchspacing \textcolor[rgb]{0.20,0.40,0.64}{\makebox[0.525em][c]{$\int$}} t\textcolor[rgb]{0.81,0.36,0.00}{\textbf{,}} \textcolor[rgb]{0.20,0.40,0.64}{...}} denotes the integral over $\R^d$ with respect to the
Lebesgue measure, and \mbox{\ttfamily\frenchspacing inner \makebox[0.525em][c]{$\mathbb{R}$} z\textcolor[rgb]{0.20,0.40,0.64}{.}\textcolor[rgb]{0,0,0.81}{\textbf{2}} t} is the inner product
$\omega\cdot t$. The annotation \mbox{\ttfamily\frenchspacing \textcolor[rgb]{0.81,0.36,0.00}{\textbf{(}}\textcolor[rgb]{0,0,0.81}{\textbf{1}} \textcolor[rgb]{0.20,0.40,0.64}{/} \textcolor[rgb]{0,0,0.81}{\textbf{2}} \textcolor[rgb]{0.81,0.36,0.00}{\textbf{:}} \makebox[0.525em][c]{$\mathbb{R}$}\textcolor[rgb]{0.81,0.36,0.00}{\textbf{)}}} tells Lean that $1/2$ is a
real number. Hence, \mbox{\ttfamily\frenchspacing ambiguity f \textcolor[rgb]{0.81,0.36,0.00}{\textbf{(}}x\textcolor[rgb]{0.81,0.36,0.00}{\textbf{,}} \makebox[0.525em][c]{$\omega$}\textcolor[rgb]{0.81,0.36,0.00}{\textbf{)}}} is precisely
$$
A(f)(x,\omega)
=\int_{\R^d} f\left(t+\frac{x}{2}\right)\overline{f\left(t-\frac{x}{2}\right)}
e^{-2\pi i\omega\cdot t}\,dt.
$$
The open zero set $\Zopen(f)$ is the interior of the
zero set of $A(f)$.

\par\medskip
{\ttfamily\frenchspacing\setlength{\parindent}{0pt}\setlength{\parskip}{0pt}\raggedright
\noindent \textcolor[rgb]{0.13,0.29,0.53}{\textbf{def}} openZeros \textcolor[rgb]{0.81,0.36,0.00}{\textbf{\{}}d \textcolor[rgb]{0.81,0.36,0.00}{\textbf{:}} \makebox[0.525em][c]{$\mathbb{N}$}\textcolor[rgb]{0.81,0.36,0.00}{\textbf{\}}} \textcolor[rgb]{0.81,0.36,0.00}{\textbf{(}}f \textcolor[rgb]{0.81,0.36,0.00}{\textbf{:}} L2 d\textcolor[rgb]{0.81,0.36,0.00}{\textbf{)}} \textcolor[rgb]{0.81,0.36,0.00}{\textbf{:}} Set \textcolor[rgb]{0.81,0.36,0.00}{\textbf{(}}PhaseSpace d\textcolor[rgb]{0.81,0.36,0.00}{\textbf{)}} \textcolor[rgb]{0.81,0.36,0.00}{\textbf{:=}}\par
\noindent \hspace*{1.05em}interior \textcolor[rgb]{0.81,0.36,0.00}{\textbf{\{}}z \textcolor[rgb]{0.20,0.40,0.64}{\textbar{}} ambiguity f z \textcolor[rgb]{0.20,0.40,0.64}{=} \textcolor[rgb]{0,0,0.81}{\textbf{0}}\textcolor[rgb]{0.81,0.36,0.00}{\textbf{\}}}\par
}\par\medskip

The interior is taken with respect to the product topology on
$\R^d\times\R^d$, which is the usual topology of $\R^{2d}$. Thus,
\mbox{\ttfamily\frenchspacing openZeros f} coincides with
$\Zopen(f)=\interior\{z\in\R^{2d}:A(f)(z)=0\}$. The short-time Fourier
transform is formalized in the same way.

\par\medskip
{\ttfamily\frenchspacing\setlength{\parindent}{0pt}\setlength{\parskip}{0pt}\raggedright
\noindent \textcolor[rgb]{0.13,0.29,0.53}{\textbf{def}} stft \textcolor[rgb]{0.81,0.36,0.00}{\textbf{\{}}d \textcolor[rgb]{0.81,0.36,0.00}{\textbf{:}} \makebox[0.525em][c]{$\mathbb{N}$}\textcolor[rgb]{0.81,0.36,0.00}{\textbf{\}}} \textcolor[rgb]{0.81,0.36,0.00}{\textbf{(}}w f \textcolor[rgb]{0.81,0.36,0.00}{\textbf{:}} L2 d\textcolor[rgb]{0.81,0.36,0.00}{\textbf{)}} \textcolor[rgb]{0.81,0.36,0.00}{\textbf{(}}z \textcolor[rgb]{0.81,0.36,0.00}{\textbf{:}} PhaseSpace d\textcolor[rgb]{0.81,0.36,0.00}{\textbf{)}} \textcolor[rgb]{0.81,0.36,0.00}{\textbf{:}} \makebox[0.525em][c]{$\mathbb{C}$} \textcolor[rgb]{0.81,0.36,0.00}{\textbf{:=}}\par
\noindent \hspace*{1.05em}\textcolor[rgb]{0.20,0.40,0.64}{\makebox[0.525em][c]{$\int$}} t\textcolor[rgb]{0.81,0.36,0.00}{\textbf{,}} f t \textcolor[rgb]{0.20,0.40,0.64}{*} star \textcolor[rgb]{0.81,0.36,0.00}{\textbf{(}}w \textcolor[rgb]{0.81,0.36,0.00}{\textbf{(}}t \textcolor[rgb]{0.20,0.40,0.64}{-} z\textcolor[rgb]{0.20,0.40,0.64}{.}\textcolor[rgb]{0,0,0.81}{\textbf{1}}\textcolor[rgb]{0.81,0.36,0.00}{\textbf{)}}\textcolor[rgb]{0.81,0.36,0.00}{\textbf{)}} \textcolor[rgb]{0.20,0.40,0.64}{*}\par
\noindent \hspace*{2.1em}Complex\textcolor[rgb]{0.20,0.40,0.64}{.}exp \textcolor[rgb]{0.81,0.36,0.00}{\textbf{(}}\textcolor[rgb]{0.20,0.40,0.64}{-}\textcolor[rgb]{0,0,0.81}{\textbf{2}} \textcolor[rgb]{0.20,0.40,0.64}{*} Real\textcolor[rgb]{0.20,0.40,0.64}{.}pi \textcolor[rgb]{0.20,0.40,0.64}{*} Complex\textcolor[rgb]{0.20,0.40,0.64}{.}I \textcolor[rgb]{0.20,0.40,0.64}{*} \textcolor[rgb]{0.81,0.36,0.00}{\textbf{(}}inner \makebox[0.525em][c]{$\mathbb{R}$} z\textcolor[rgb]{0.20,0.40,0.64}{.}\textcolor[rgb]{0,0,0.81}{\textbf{2}} t \textcolor[rgb]{0.81,0.36,0.00}{\textbf{:}} \makebox[0.525em][c]{$\mathbb{R}$}\textcolor[rgb]{0.81,0.36,0.00}{\textbf{)}}\textcolor[rgb]{0.81,0.36,0.00}{\textbf{)}}\par
}\par\medskip

Note that the window is the first argument of \mbox{\ttfamily\frenchspacing stft}, so that
\mbox{\ttfamily\frenchspacing stft w f \textcolor[rgb]{0.81,0.36,0.00}{\textbf{(}}x\textcolor[rgb]{0.81,0.36,0.00}{\textbf{,}} \makebox[0.525em][c]{$\omega$}\textcolor[rgb]{0.81,0.36,0.00}{\textbf{)}}} equals $V_wf(x,\omega)$. The phase retrieval property
of a window is defined in terms of \mbox{\ttfamily\frenchspacing stft}.

\par\medskip
{\ttfamily\frenchspacing\setlength{\parindent}{0pt}\setlength{\parskip}{0pt}\raggedright
\noindent \textcolor[rgb]{0.13,0.29,0.53}{\textbf{def}} DoesPhaseRetrieval \textcolor[rgb]{0.81,0.36,0.00}{\textbf{\{}}d \textcolor[rgb]{0.81,0.36,0.00}{\textbf{:}} \makebox[0.525em][c]{$\mathbb{N}$}\textcolor[rgb]{0.81,0.36,0.00}{\textbf{\}}} \textcolor[rgb]{0.81,0.36,0.00}{\textbf{(}}w \textcolor[rgb]{0.81,0.36,0.00}{\textbf{:}} L2 d\textcolor[rgb]{0.81,0.36,0.00}{\textbf{)}} \textcolor[rgb]{0.81,0.36,0.00}{\textbf{:}} \textcolor[rgb]{0.13,0.29,0.53}{\textbf{Prop}} \textcolor[rgb]{0.81,0.36,0.00}{\textbf{:=}}\par
\noindent \hspace*{1.05em}\textcolor[rgb]{0.20,0.40,0.64}{\makebox[0.525em][c]{$\forall$}} f g \textcolor[rgb]{0.81,0.36,0.00}{\textbf{:}} L2 d\textcolor[rgb]{0.81,0.36,0.00}{\textbf{,}}\par
\noindent \hspace*{2.1em}\textcolor[rgb]{0.81,0.36,0.00}{\textbf{(}}\textcolor[rgb]{0.20,0.40,0.64}{\makebox[0.525em][c]{$\forall$}} z\textcolor[rgb]{0.81,0.36,0.00}{\textbf{,}} \textcolor[rgb]{0.20,0.40,0.64}{\makebox[0.525em][c]{$\Vert$}}stft w f z\textcolor[rgb]{0.20,0.40,0.64}{\makebox[0.525em][c]{$\Vert$}} \textcolor[rgb]{0.20,0.40,0.64}{=} \textcolor[rgb]{0.20,0.40,0.64}{\makebox[0.525em][c]{$\Vert$}}stft w g z\textcolor[rgb]{0.20,0.40,0.64}{\makebox[0.525em][c]{$\Vert$}}\textcolor[rgb]{0.81,0.36,0.00}{\textbf{)}} \textcolor[rgb]{0.20,0.40,0.64}{\makebox[0.525em][c]{$\rightarrow$}}\par
\noindent \hspace*{2.1em}\textcolor[rgb]{0.20,0.40,0.64}{\makebox[0.525em][c]{$\exists$}} \makebox[0.525em][c]{$\theta$} \textcolor[rgb]{0.81,0.36,0.00}{\textbf{:}} \makebox[0.525em][c]{$\mathbb{R}$}\textcolor[rgb]{0.81,0.36,0.00}{\textbf{,}} g \textcolor[rgb]{0.20,0.40,0.64}{=} Complex\textcolor[rgb]{0.20,0.40,0.64}{.}exp \textcolor[rgb]{0.81,0.36,0.00}{\textbf{(}}\makebox[0.525em][c]{$\theta$} \textcolor[rgb]{0.20,0.40,0.64}{*} Complex\textcolor[rgb]{0.20,0.40,0.64}{.}I\textcolor[rgb]{0.81,0.36,0.00}{\textbf{)}} \textcolor[rgb]{0.20,0.40,0.64}{\makebox[0.525em][c]{$\bullet$}} f\par
}\par\medskip

Since \mbox{\ttfamily\frenchspacing \textcolor[rgb]{0.20,0.40,0.64}{\makebox[0.525em][c]{$\Vert$}\makebox[0.525em][c]{$\cdot$}\makebox[0.525em][c]{$\Vert$}}} is the modulus of a complex number, the predicate
\mbox{\ttfamily\frenchspacing DoesPhaseRetrieval w} states that for all $f,g\in L^2(\R^d)$ the
identity $|V_wf|=|V_wg|$ on $\R^{2d}$ implies $g=e^{i\theta}f$ for some
$\theta\in\R$. This is exactly the statement that $w$ does phase retrieval.

\subsection{Formalization of Theorem~\ref{thm:phase-retrieval-d}}

Theorem~\ref{thm:phase-retrieval-d} is stated in Lean as follows.

\par\medskip
{\ttfamily\frenchspacing\setlength{\parindent}{0pt}\setlength{\parskip}{0pt}\raggedright
\noindent \textcolor[rgb]{0.13,0.29,0.53}{\textbf{theorem}} bounded\_openZeros\_phaseRetrieval\par
\noindent \hspace*{2.1em}\textcolor[rgb]{0.81,0.36,0.00}{\textbf{\{}}d \textcolor[rgb]{0.81,0.36,0.00}{\textbf{:}} \makebox[0.525em][c]{$\mathbb{N}$}\textcolor[rgb]{0.81,0.36,0.00}{\textbf{\}}} \textcolor[rgb]{0.81,0.36,0.00}{\textbf{(}}hd \textcolor[rgb]{0.81,0.36,0.00}{\textbf{:}} \textcolor[rgb]{0,0,0.81}{\textbf{1}} \textcolor[rgb]{0.20,0.40,0.64}{\makebox[0.525em][c]{$\leq$}} d\textcolor[rgb]{0.81,0.36,0.00}{\textbf{)}} \textcolor[rgb]{0.81,0.36,0.00}{\textbf{(}}w \textcolor[rgb]{0.81,0.36,0.00}{\textbf{:}} L2 d\textcolor[rgb]{0.81,0.36,0.00}{\textbf{)}} \textcolor[rgb]{0.81,0.36,0.00}{\textbf{(}}hw \textcolor[rgb]{0.81,0.36,0.00}{\textbf{:}} w \textcolor[rgb]{0.20,0.40,0.64}{\makebox[0.525em][c]{$\neq$}} \textcolor[rgb]{0,0,0.81}{\textbf{0}}\textcolor[rgb]{0.81,0.36,0.00}{\textbf{)}}\par
\noindent \hspace*{2.1em}\textcolor[rgb]{0.81,0.36,0.00}{\textbf{(}}hbounded \textcolor[rgb]{0.81,0.36,0.00}{\textbf{:}} Bornology\textcolor[rgb]{0.20,0.40,0.64}{.}IsBounded \textcolor[rgb]{0.81,0.36,0.00}{\textbf{(}}openZeros w\textcolor[rgb]{0.81,0.36,0.00}{\textbf{)}}\textcolor[rgb]{0.81,0.36,0.00}{\textbf{)}} \textcolor[rgb]{0.81,0.36,0.00}{\textbf{:}}\par
\noindent \hspace*{2.1em}DoesPhaseRetrieval w\par
}\par\medskip

The dimension $d\geq1$ is an implicit natural number, the window $w$ is a
nonzero element of $L^2(\R^d)$, and \mbox{\ttfamily\frenchspacing Bornology\textcolor[rgb]{0.20,0.40,0.64}{.}IsBounded} is Mathlib's
notion of a bounded set. Thus, the theorem states that if
$w\in L^2(\R^d)\setminus\{0\}$ and $\Zopen(w)$ is bounded, then $w$ does
phase retrieval.

\subsection{Formalization of Theorem~\ref{thm:main}}

Theorem~\ref{thm:main} is stated in Lean as follows.

\par\medskip
{\ttfamily\frenchspacing\setlength{\parindent}{0pt}\setlength{\parskip}{0pt}\raggedright
\noindent \textcolor[rgb]{0.13,0.29,0.53}{\textbf{theorem}} prescribed\_openZeros\par
\noindent \hspace*{2.1em}\textcolor[rgb]{0.81,0.36,0.00}{\textbf{\{}}d \textcolor[rgb]{0.81,0.36,0.00}{\textbf{:}} \makebox[0.525em][c]{$\mathbb{N}$}\textcolor[rgb]{0.81,0.36,0.00}{\textbf{\}}} \textcolor[rgb]{0.81,0.36,0.00}{\textbf{(}}hd \textcolor[rgb]{0.81,0.36,0.00}{\textbf{:}} \textcolor[rgb]{0,0,0.81}{\textbf{1}} \textcolor[rgb]{0.20,0.40,0.64}{\makebox[0.525em][c]{$\leq$}} d\textcolor[rgb]{0.81,0.36,0.00}{\textbf{)}} \textcolor[rgb]{0.81,0.36,0.00}{\textbf{(}}\makebox[0.525em][c]{$\Omega$} \textcolor[rgb]{0.81,0.36,0.00}{\textbf{:}} Set \textcolor[rgb]{0.81,0.36,0.00}{\textbf{(}}PhaseSpace d\textcolor[rgb]{0.81,0.36,0.00}{\textbf{)}}\textcolor[rgb]{0.81,0.36,0.00}{\textbf{)}}\par
\noindent \hspace*{2.1em}\textcolor[rgb]{0.81,0.36,0.00}{\textbf{(}}hbounded \textcolor[rgb]{0.81,0.36,0.00}{\textbf{:}} Bornology\textcolor[rgb]{0.20,0.40,0.64}{.}IsBounded \makebox[0.525em][c]{$\Omega$}\textcolor[rgb]{0.81,0.36,0.00}{\textbf{)}} \textcolor[rgb]{0.81,0.36,0.00}{\textbf{(}}hne \textcolor[rgb]{0.81,0.36,0.00}{\textbf{:}} \makebox[0.525em][c]{$\Omega$}\textcolor[rgb]{0.20,0.40,0.64}{.}Nonempty\textcolor[rgb]{0.81,0.36,0.00}{\textbf{)}}\par
\noindent \hspace*{2.1em}\textcolor[rgb]{0.81,0.36,0.00}{\textbf{(}}hregular \textcolor[rgb]{0.81,0.36,0.00}{\textbf{:}} \makebox[0.525em][c]{$\Omega$} \textcolor[rgb]{0.20,0.40,0.64}{=} interior \textcolor[rgb]{0.81,0.36,0.00}{\textbf{(}}closure \makebox[0.525em][c]{$\Omega$}\textcolor[rgb]{0.81,0.36,0.00}{\textbf{)}}\textcolor[rgb]{0.81,0.36,0.00}{\textbf{)}}\par
\noindent \hspace*{2.1em}\textcolor[rgb]{0.81,0.36,0.00}{\textbf{(}}a \textcolor[rgb]{0.81,0.36,0.00}{\textbf{:}} PhaseSpace d\textcolor[rgb]{0.81,0.36,0.00}{\textbf{)}} \textcolor[rgb]{0.81,0.36,0.00}{\textbf{(}}ha \textcolor[rgb]{0.81,0.36,0.00}{\textbf{:}} a \textcolor[rgb]{0.20,0.40,0.64}{\makebox[0.525em][c]{$\neq$}} \textcolor[rgb]{0,0,0.81}{\textbf{0}}\textcolor[rgb]{0.81,0.36,0.00}{\textbf{)}} \textcolor[rgb]{0.81,0.36,0.00}{\textbf{:}}\par
\noindent \hspace*{2.1em}\textcolor[rgb]{0.20,0.40,0.64}{\makebox[0.525em][c]{$\exists$}} b\makebox[0.525em][c]{${}_0$} \textcolor[rgb]{0.20,0.40,0.64}{\textgreater{}} \textcolor[rgb]{0.81,0.36,0.00}{\textbf{(}}\textcolor[rgb]{0,0,0.81}{\textbf{0}} \textcolor[rgb]{0.81,0.36,0.00}{\textbf{:}} \makebox[0.525em][c]{$\mathbb{R}$}\textcolor[rgb]{0.81,0.36,0.00}{\textbf{)}}\textcolor[rgb]{0.81,0.36,0.00}{\textbf{,}} \textcolor[rgb]{0.20,0.40,0.64}{\makebox[0.525em][c]{$\forall$}} b \textcolor[rgb]{0.81,0.36,0.00}{\textbf{:}} \makebox[0.525em][c]{$\mathbb{R}$}\textcolor[rgb]{0.81,0.36,0.00}{\textbf{,}} \textcolor[rgb]{0,0,0.81}{\textbf{0}} \textcolor[rgb]{0.20,0.40,0.64}{\textless{}} b \textcolor[rgb]{0.20,0.40,0.64}{\makebox[0.525em][c]{$\rightarrow$}} b \textcolor[rgb]{0.20,0.40,0.64}{\textless{}} b\makebox[0.525em][c]{${}_0$} \textcolor[rgb]{0.20,0.40,0.64}{\makebox[0.525em][c]{$\rightarrow$}}\par
\noindent \hspace*{3.1500000000000004em}\textcolor[rgb]{0.20,0.40,0.64}{\makebox[0.525em][c]{$\exists$}} f \textcolor[rgb]{0.81,0.36,0.00}{\textbf{:}} L2 d\textcolor[rgb]{0.81,0.36,0.00}{\textbf{,}}\par
\noindent \hspace*{4.2em}openZeros f \textcolor[rgb]{0.20,0.40,0.64}{=}\par
\noindent \hspace*{5.25em}\textcolor[rgb]{0.81,0.36,0.00}{\textbf{(}}\textcolor[rgb]{0.13,0.29,0.53}{\textbf{fun}} z \textcolor[rgb]{0.20,0.40,0.64}{\makebox[0.525em][c]{$\mapsto$}} a \textcolor[rgb]{0.20,0.40,0.64}{+} b \textcolor[rgb]{0.20,0.40,0.64}{\makebox[0.525em][c]{$\bullet$}} z\textcolor[rgb]{0.81,0.36,0.00}{\textbf{)}} \textcolor[rgb]{0.20,0.40,0.64}{\textquotesingle{}}\textcolor[rgb]{0.20,0.40,0.64}{\textquotesingle{}} \makebox[0.525em][c]{$\Omega$} \textcolor[rgb]{0.20,0.40,0.64}{\makebox[0.525em][c]{$\cup$}} \textcolor[rgb]{0.81,0.36,0.00}{\textbf{(}}\textcolor[rgb]{0.13,0.29,0.53}{\textbf{fun}} z \textcolor[rgb]{0.20,0.40,0.64}{\makebox[0.525em][c]{$\mapsto$}} \textcolor[rgb]{0.20,0.40,0.64}{-}\textcolor[rgb]{0.81,0.36,0.00}{\textbf{(}}a \textcolor[rgb]{0.20,0.40,0.64}{+} b \textcolor[rgb]{0.20,0.40,0.64}{\makebox[0.525em][c]{$\bullet$}} z\textcolor[rgb]{0.81,0.36,0.00}{\textbf{)}}\textcolor[rgb]{0.81,0.36,0.00}{\textbf{)}} \textcolor[rgb]{0.20,0.40,0.64}{\textquotesingle{}}\textcolor[rgb]{0.20,0.40,0.64}{\textquotesingle{}} \makebox[0.525em][c]{$\Omega$}\par
}\par\medskip

The hypotheses \mbox{\ttfamily\frenchspacing hbounded} and \mbox{\ttfamily\frenchspacing hne} state that $\Omega$ is bounded
and nonempty. The hypothesis \mbox{\ttfamily\frenchspacing hregular} states that $\Omega$ is regular
open, i.e., $\Omega$ coincides with the interior of its closure. The center
$a$ is a nonzero point of phase space. The expression
\mbox{\ttfamily\frenchspacing \textcolor[rgb]{0.20,0.40,0.64}{\makebox[0.525em][c]{$\exists$}} b\makebox[0.525em][c]{${}_0$} \textcolor[rgb]{0.20,0.40,0.64}{\textgreater{}} \textcolor[rgb]{0.81,0.36,0.00}{\textbf{(}}\textcolor[rgb]{0,0,0.81}{\textbf{0}} \textcolor[rgb]{0.81,0.36,0.00}{\textbf{:}} \makebox[0.525em][c]{$\mathbb{R}$}\textcolor[rgb]{0.81,0.36,0.00}{\textbf{),}} P} abbreviates the existence of a real number $b_0>0$
with property $P$. Further, \mbox{\ttfamily\frenchspacing \textcolor[rgb]{0.81,0.36,0.00}{\textbf{(}}\textcolor[rgb]{0.13,0.29,0.53}{\textbf{fun}} z \textcolor[rgb]{0.20,0.40,0.64}{\makebox[0.525em][c]{$\mapsto$}} a \textcolor[rgb]{0.20,0.40,0.64}{+} b \textcolor[rgb]{0.20,0.40,0.64}{\makebox[0.525em][c]{$\bullet$}} z\textcolor[rgb]{0.81,0.36,0.00}{\textbf{)}} \textcolor[rgb]{0.20,0.40,0.64}{\textquotesingle{}\textquotesingle{}} \makebox[0.525em][c]{$\Omega$}} is the image of
$\Omega$ under the map $z\mapsto a+bz$, that is, the set $a+b\Omega$, and
similarly \mbox{\ttfamily\frenchspacing \textcolor[rgb]{0.81,0.36,0.00}{\textbf{(}}\textcolor[rgb]{0.13,0.29,0.53}{\textbf{fun}} z \textcolor[rgb]{0.20,0.40,0.64}{\makebox[0.525em][c]{$\mapsto$}} \textcolor[rgb]{0.20,0.40,0.64}{-}\textcolor[rgb]{0.81,0.36,0.00}{\textbf{(}}a \textcolor[rgb]{0.20,0.40,0.64}{+} b \textcolor[rgb]{0.20,0.40,0.64}{\makebox[0.525em][c]{$\bullet$}} z\textcolor[rgb]{0.81,0.36,0.00}{\textbf{))}} \textcolor[rgb]{0.20,0.40,0.64}{\textquotesingle{}\textquotesingle{}} \makebox[0.525em][c]{$\Omega$}} is the set $-(a+b\Omega)$. Since
$b_0$ is chosen after $d$, $\Omega$ and $a$ are fixed, it depends only on these
data, as in Theorem~\ref{thm:main}. Hence, the conclusion states that there
exists $b_0>0$ such that for every $0<b<b_0$ one can find $f\in L^2(\R^d)$
with
$$
\Zopen(f)=(a+b\Omega)\cup(-(a+b\Omega)).
$$

\subsection{Formalization of Theorem~\ref{thm:counterexample}}

The annulus $O_{r,R}$ defined in \eqref{eq:annulus} is formalized as follows.

\par\medskip
{\ttfamily\frenchspacing\setlength{\parindent}{0pt}\setlength{\parskip}{0pt}\raggedright
\noindent \textcolor[rgb]{0.13,0.29,0.53}{\textbf{def}} annulus \textcolor[rgb]{0.81,0.36,0.00}{\textbf{(}}d \textcolor[rgb]{0.81,0.36,0.00}{\textbf{:}} \makebox[0.525em][c]{$\mathbb{N}$}\textcolor[rgb]{0.81,0.36,0.00}{\textbf{)}} \textcolor[rgb]{0.81,0.36,0.00}{\textbf{(}}r R \textcolor[rgb]{0.81,0.36,0.00}{\textbf{:}} \makebox[0.525em][c]{$\mathbb{R}$}\textcolor[rgb]{0.81,0.36,0.00}{\textbf{)}} \textcolor[rgb]{0.81,0.36,0.00}{\textbf{:}} Set \textcolor[rgb]{0.81,0.36,0.00}{\textbf{(}}PhaseSpace d\textcolor[rgb]{0.81,0.36,0.00}{\textbf{)}} \textcolor[rgb]{0.81,0.36,0.00}{\textbf{:=}}\par
\noindent \hspace*{1.05em}\textcolor[rgb]{0.81,0.36,0.00}{\textbf{\{}}z \textcolor[rgb]{0.20,0.40,0.64}{\textbar{}} r \textcolor[rgb]{0.20,0.40,0.64}{\textless{}} Real\textcolor[rgb]{0.20,0.40,0.64}{.}sqrt \textcolor[rgb]{0.81,0.36,0.00}{\textbf{(}}\textcolor[rgb]{0.20,0.40,0.64}{\makebox[0.525em][c]{$\Vert$}}z\textcolor[rgb]{0.20,0.40,0.64}{.}\textcolor[rgb]{0,0,0.81}{\textbf{1}}\textcolor[rgb]{0.20,0.40,0.64}{\makebox[0.525em][c]{$\Vert$}} \textcolor[rgb]{0.20,0.40,0.64}{\textasciicircum{}} \textcolor[rgb]{0,0,0.81}{\textbf{2}} \textcolor[rgb]{0.20,0.40,0.64}{+} \textcolor[rgb]{0.20,0.40,0.64}{\makebox[0.525em][c]{$\Vert$}}z\textcolor[rgb]{0.20,0.40,0.64}{.}\textcolor[rgb]{0,0,0.81}{\textbf{2}}\textcolor[rgb]{0.20,0.40,0.64}{\makebox[0.525em][c]{$\Vert$}} \textcolor[rgb]{0.20,0.40,0.64}{\textasciicircum{}} \textcolor[rgb]{0,0,0.81}{\textbf{2}}\textcolor[rgb]{0.81,0.36,0.00}{\textbf{)}} \textcolor[rgb]{0.20,0.40,0.64}{\makebox[0.525em][c]{$\wedge$}}\par
\noindent \hspace*{3.6750000000000003em}Real\textcolor[rgb]{0.20,0.40,0.64}{.}sqrt \textcolor[rgb]{0.81,0.36,0.00}{\textbf{(}}\textcolor[rgb]{0.20,0.40,0.64}{\makebox[0.525em][c]{$\Vert$}}z\textcolor[rgb]{0.20,0.40,0.64}{.}\textcolor[rgb]{0,0,0.81}{\textbf{1}}\textcolor[rgb]{0.20,0.40,0.64}{\makebox[0.525em][c]{$\Vert$}} \textcolor[rgb]{0.20,0.40,0.64}{\textasciicircum{}} \textcolor[rgb]{0,0,0.81}{\textbf{2}} \textcolor[rgb]{0.20,0.40,0.64}{+} \textcolor[rgb]{0.20,0.40,0.64}{\makebox[0.525em][c]{$\Vert$}}z\textcolor[rgb]{0.20,0.40,0.64}{.}\textcolor[rgb]{0,0,0.81}{\textbf{2}}\textcolor[rgb]{0.20,0.40,0.64}{\makebox[0.525em][c]{$\Vert$}} \textcolor[rgb]{0.20,0.40,0.64}{\textasciicircum{}} \textcolor[rgb]{0,0,0.81}{\textbf{2}}\textcolor[rgb]{0.81,0.36,0.00}{\textbf{)}} \textcolor[rgb]{0.20,0.40,0.64}{\textless{}} R\textcolor[rgb]{0.81,0.36,0.00}{\textbf{\}}}\par
}\par\medskip

Mathlib equips the product $\R^d\times\R^d$ with the maximum norm
$\max\{|x|,|\omega|\}$. For this reason, the Euclidean norm
$|z|=(|x|^2+|\omega|^2)^{1/2}$ of $z=(x,\omega)$ is written out explicitly,
and \mbox{\ttfamily\frenchspacing annulus d r R} coincides with $O_{r,R}$.
Theorem~\ref{thm:counterexample} is stated in Lean as follows.

\par\medskip
{\ttfamily\frenchspacing\setlength{\parindent}{0pt}\setlength{\parskip}{0pt}\raggedright
\noindent \textcolor[rgb]{0.13,0.29,0.53}{\textbf{theorem}} forbidden\_annulus \textcolor[rgb]{0.81,0.36,0.00}{\textbf{:}}\par
\noindent \hspace*{2.1em}\textcolor[rgb]{0.20,0.40,0.64}{\makebox[0.525em][c]{$\exists$}} r R \textcolor[rgb]{0.81,0.36,0.00}{\textbf{:}} \makebox[0.525em][c]{$\mathbb{R}$}\textcolor[rgb]{0.81,0.36,0.00}{\textbf{,}} \textcolor[rgb]{0,0,0.81}{\textbf{0}} \textcolor[rgb]{0.20,0.40,0.64}{\textless{}} r \textcolor[rgb]{0.20,0.40,0.64}{\makebox[0.525em][c]{$\wedge$}} r \textcolor[rgb]{0.20,0.40,0.64}{\textless{}} R \textcolor[rgb]{0.20,0.40,0.64}{\makebox[0.525em][c]{$\wedge$}}\par
\noindent \hspace*{3.1500000000000004em}\textcolor[rgb]{0.20,0.40,0.64}{\makebox[0.525em][c]{$\forall$}} f \textcolor[rgb]{0.81,0.36,0.00}{\textbf{:}} L2 \textcolor[rgb]{0,0,0.81}{\textbf{1}}\textcolor[rgb]{0.81,0.36,0.00}{\textbf{,}} f \textcolor[rgb]{0.20,0.40,0.64}{\makebox[0.525em][c]{$\neq$}} \textcolor[rgb]{0,0,0.81}{\textbf{0}} \textcolor[rgb]{0.20,0.40,0.64}{\makebox[0.525em][c]{$\rightarrow$}} \textcolor[rgb]{0.20,0.40,0.64}{\makebox[0.525em][c]{$\neg$}} annulus \textcolor[rgb]{0,0,0.81}{\textbf{1}} r R \textcolor[rgb]{0.20,0.40,0.64}{\makebox[0.525em][c]{$\subseteq$}} openZeros f\par
}\par\medskip

Here \mbox{\ttfamily\frenchspacing L2 \textcolor[rgb]{0,0,0.81}{\textbf{1}}} is the space $L^2(\R)$, and
\mbox{\ttfamily\frenchspacing \textcolor[rgb]{0.20,0.40,0.64}{\makebox[0.525em][c]{$\neg$}} annulus \textcolor[rgb]{0,0,0.81}{\textbf{1}} r R \textcolor[rgb]{0.20,0.40,0.64}{\makebox[0.525em][c]{$\subseteq$}} openZeros f} states that $O_{r,R}$ is not contained in
$\Zopen(f)$.

\subsection{Source code}

The repository separates the statements from the proofs, so that the
formalized claims can be inspected independently. The file
\mbox{\ttfamily\frenchspacing Showcase\textcolor[rgb]{0.20,0.40,0.64}{.}lean} contains the definitions and statements described above,
with the three proofs replaced by the placeholder \mbox{\ttfamily\frenchspacing \textcolor[rgb]{0.94,0.16,0.16}{sorry}}. Verifying that
the formalization faithfully captures
Theorems~\ref{thm:phase-retrieval-d}--\ref{thm:counterexample} therefore only
requires reading this file. The companion file \mbox{\ttfamily\frenchspacing Showcase\_WithProofs\textcolor[rgb]{0.20,0.40,0.64}{.}lean}
contains the same statements, but each placeholder is replaced by a complete
proof, whose correctness is certified by the Lean kernel during compilation. A
separate script checks that the proved theorems have exactly the statements of
\mbox{\ttfamily\frenchspacing Showcase\textcolor[rgb]{0.20,0.40,0.64}{.}lean}. Moreover, Lean's axiom audit shows that the three proofs
depend only on the standard axioms \mbox{\ttfamily\frenchspacing propext}, \mbox{\ttfamily\frenchspacing Classical\textcolor[rgb]{0.20,0.40,0.64}{.}choice} and
\mbox{\ttfamily\frenchspacing Quot\textcolor[rgb]{0.20,0.40,0.64}{.}sound}.

\section*{Declaration of AI-assisted tools}

During the preparation of this manuscript, the authors made use of generative AI tools for a variety of purposes, including improving the presentation, exploring mathematical ideas, suggesting examples and possible approaches to proofs, assisting with computations, and preparing illustrative material. All AI-generated mathematical content was critically assessed, and the final statements, proofs, and computations were independently checked and rewritten by the authors. More specifically, the existence of an ambiguity function vanishing on a bounded open set was first produced by GPT-5.4. Discussions with GPT-5.5 and 5.6 were important for obtaining the final, general results in this manuscript. The initial Lean verification was done by Opus 4.6 and newer models were only used for improving the readability of the paper and for additional Lean verifications.
The authors take full responsibility for the accuracy, originality, and integrity of the manuscript.

\section*{Acknowledgments}

P.~J.~is partially supported by  the French National Research Agency (ANR) under contract number
ANR-24-CE40-5470. This project has received financial support from the French government within the framework of the France 2030 programme IDEX Université de Bordeaux.

L.~L.~is grateful to the Azrieli Foundation for the award of an Azrieli Fellowship and acknowledges the support of this research by ISF Grant No.~854/25.

We would like to thank Rima Alaifari and Victoire Laurent for interesting discussions on ambiguity functions.

\bibliographystyle{plain}

\end{document}